\documentclass[12pt]{article}
\usepackage{amsmath,amssymb}
\usepackage{graphicx}
\newcommand{\eqdef}{\stackrel{\text{def}}{=}}

\newcommand{\n}{\nonumber \\}
\newcommand{\bm}{\boldsymbol}
\newcommand{\ignore}[1]{}
\numberwithin{equation}{section}
\newcommand{\Romannumeral}[1]{\uppercase\expandafter{\romannumeral#1}}

\newtheorem{theo}{\bf Theorem}[section]

\newtheorem{rema}[theo]{\bf Remark}

\newtheorem{prop}[theo]{\bf Proposition}

\allowdisplaybreaks[4]

\begin{document}

\baselineskip=20pt
%%%%%%%%%%%%%%%%%%%%%%%%%%%%%%%%%%%%%%%%%%%%%%%%%%%%%%%%%%%%
%                                                          %
%  Title page                                              %
%                                                          %
%%%%%%%%%%%%%%%%%%%%%%%%%%%%%%%%%%%%%%%%%%%%%%%%%%%%%%%%%%%%
\newcommand{\preprint}{
\vspace*{-20mm}\begin{flushleft}\end{flushleft}
}
\newcommand{\Title}[1]{{\baselineskip=26pt
  \begin{center} \Large \bf #1 \\ \ \\ \end{center}}}
\newcommand{\Author}{\begin{center}
  \large \bf 
  Ryu Sasaki${}$ \end{center}}
\newcommand{\Address}{\begin{center}
     Department of Physics and Astronomy, Tokyo University of Science,
     Noda 278-8510, Japan
        \end{center}}
\newcommand{\Accepted}[1]{\begin{center}
  {\large \sf #1}\\ \vspace{1mm}{\small \sf Accepted for Publication}
  \end{center}}

\preprint
\thispagestyle{empty}

\Title{Multivariate Kawtchouk polynomials as  Birth and Death polynomials
}

\Author

\Address
\vspace{1cm}

\begin{abstract}
Multivariate Krawtchouk polynomials are constructed explicitly as Birth and Death polynomials,
which have the nearest neighbour interactions.
They form the complete set of eigenpolynomials of a birth and death process 
with the birth and death rates at population 
$\bm{x}=(x_1,\ldots,x_n)$ are $B_j(\bm{x})=\bigl(N-\sum_{i=1}^nx_i\bigr)$ 
and $D_j(\bm{x})=p_i^{-1}x_j$,  $0<p_j$, $j=1,\ldots,n$. The corresponding 
stationary distribution is the multinomial distribution with the probabilities $\{\eta_i\}$, 
$\eta_i\eqdef p_i/(1+\sum_{j=1}^np_j)$.
The polynomials, depending on $n+1$ parameters ($\{p_i\}$ and $N$),
 satisfy the difference equation with the coefficients 
$B_j(\bm{x})$ and $D_j(\bm{x})$  $j=1,\ldots,n$,
which is the straightforward generalisation of the difference equation 
governing the single variable Krawtchouk polynomials.
The polynomials  are  truncated $(n+1,2n+2)$ hypergeometric functions of Aomoto-Gelfand.
The divariate Rahman polynomials are identified as the dual polynomials with a special parametrisation.
\end{abstract}

%%%%%%%%%%%%%%%%%%%%%%%%%%%%%%%%%%%%%%%%%%%%
%
%   Section1 Introduction
%
%%%%%%%%%%%%%%%%%%%%%%%%%%%%%%%%%%%%%%%%%%%%
\section{Introduction}
\label{sec:intro}
As a first member of the discrete multivariate hypergeometric orthogonal polynomials of Askey scheme
\cite{askey,ismail, koeswart,gasper}, satisfying second order difference equations with nearest neighbour
interactions, the multivariate Krawtchouk polynomials are constructed {\em explicitly} as
the simplest multivariate Birth and Death (BD) \cite{feller, KarMcG, ismail} polynomials.
Multivariate BD processes are the nearest neighbour interactions of multi-dimensional discrete systems,
the most basic type of interactions, like the well known Ising models.

Compared with the single variable cases, construction of multivariate orthogonal polynomials
has much varieties.  Roughly speaking, single variable orthogonal polynomials are almost uniquely determined
when the orthogonality measure is given.
In contrast, formulations of multivariate orthogonal polynomials require additional inputs 
other than the orthogonality measures.
Products of single variable polynomials are used in many examples \cite{tra, illxu, khare, xu}.
The well known Rahman polynomials are generated as the eigenfunctions of certain Markov processes
related with cumulative Bernoulli trials having the multinomial distributions \cite{coo-hoa-rah77,HR, gr1, gr2, gr3, IT, Gri2}.
Connection with other orthogonal objects, {\em e.g.} 9-$j$ symbols, some group or algebra representations
are also utilised \cite{zheda, IT, I, mizu1, mt, genest}. 
Starting from some general hypergeometric functions \cite{AK, gelfand} 
 and imposing conditions related with the orthogonality 
measures is another path \cite{Gri1,Gri2,diaconis13}.

The path adopted by this paper is different from these.
Each member of the Askey scheme of discrete orthogonal polynomials satisfies a difference equation 
 with the coefficient function $B(x)$ and $D(x)$ \eqref{bdeq0} on top of the three term 
 recurrence relations.
 By a similarity transformation in terms of the corresponding orthogonality weight $W(x)$ \eqref{BDW},
 the difference equation is mapped to an operator governing the  birth and death process \eqref{BDeq}
 with the birth rate at population $x$ is $B(x)$ and death rate $D(x)$, which is exactly solvable by
 construction \cite{bdsol}.
 By reversing the logic and applying it to multivariate cases,
 one could hope that a good birth and death process for $n$ population groups with $n$ birth rates
 $\{B_i(\bm{x})\}$ and death rates  $\{D_i(\bm{x})\}$ would determine a stationary distribution
 $W(\bm{x})$ which leads to a {\em solvable} difference equation for the eigenpolynomials 
 by the inverse similarity transformation.
 However, the consistency of the interactions in various directions for obtaining the stationary 
 distribution and the solvability 
 of the difference equation put severe restrictions on the birth/death rates.
 Here I present the first successful example, multivariate Krawtchouk polynomials.
 
 This paper is organised as follows. In section \ref{sec:onedim}
the basic relation between the birth and death equation and the corresponding 
difference equation in the single variable case is recapitulated as {\bf Theorem \ref{theo1}}.
In section \ref{sec:genset}, starting with the general formulation of multivariate BD processes,
it is pointed out in {\bf Proposition \ref{compcond}} that the existence of the stationary distribution imposes severe conditions on the birth/death rates. {\bf Proposition \ref{semidef}} states that the similarity transformation of the
BD operator $L_{BD}$ \eqref{LBDop0} in terms of the square root of the stationary distribution \eqref{Hdef}
produces a real symmetric and positive semidefinte matrix $\mathcal{H}$ \eqref{Hfac}.
The multivariate counterpart of the difference operator $\widetilde{\mathcal H}$ \eqref{Hthdef2}
is introduced by another similarity transformation by $\sqrt{W(\bm{x})}$.
After reviewing the properties of the single variable Krawtchouk polynomials in section 3, 
the main results of the multivariate Krawtchouk polynomials are stated in {\bf Theorem \ref{theo:main}}.
The polynomial \eqref{Pm} is the truncated Aomoto-Gelfand \cite{AK, gelfand} hypergeometric function 
of type $(n+1,2n+2)$ which is a natural generalisation 
of the single variable formula \eqref{Krform}.
The detailed derivation follows.
In section four, the bivariate Rahman polynomials \cite{gr2} are rederived as the dual birth and death polynomials
with a very special parametrisation of the birth and death rates \eqref{dpardef1}, \eqref{dpardef2}.
The final section is for some comments. An expert can go directly to section \ref{sec:Main}.

%%%%%%%%%%%%%%%%%%%%%%%%%%%%%%%%%%%%%%%%%%%%
%
%   Section2 Birth and Death processes
%
%%%%%%%%%%%%%%%%%%%%%%%%%%%%%%%%%%%%%%%%%%%%
\section{Birth and Death processes}
\label{sec:BD}

%%%%%%%%%%%%%%%%%%%%%%%%%%%%%%%%%%%%%%%%%%%%%%%%%%%%%
%                                                   %
% 2.1 Single variable case                               %
%                                                   %
%%%%%%%%%%%%%%%%%%%%%%%%%%%%%%%%%%%%%%%%%%%%%%%%%%%%%
\subsection{Single variable case}
\label{sec:onedim}

Let us start with the review of the one-dimensional Birth and Death (BD) processes \cite{feller,ismail}.
Let $\mathcal{P}(x;t)$ be the probability distribution at time $t$ over a population space $\mathcal{X}$,
which is  a one-dimensional  nonnegative integer lattice, either finite or semi-infinite:
\begin{equation*}
\mathcal{X}=\{0,1,\ldots,N\}: \quad \text{finite},\qquad \mathcal{X}
={\mathbb N_0}: \quad \text {semi-infinite}.
\end{equation*}
Here ${\mathbb N_0}$ is the set of nonnegative integers. 
Let us denote the birth rate at population $x$ by $B(x)>0$ and the death rate by $D(x)>0$.
The time evolution of the probability distribution is governed by the following differential equation:
\begin{align}
\frac{\partial}{\partial t}\mathcal{P}(x;t)&=(L_{BD}\mathcal{P})(x;t)
=\sum_{y\in\mathcal{X}}{L_{BD}}_{x\,y}\mathcal{P}(y;t),
\quad \mathcal{P}(x;t)\ge0,\quad \sum_{x\in\mathcal{X}} \mathcal{P}(x;t)=1,
\label{bdeqformal}\\
&=-\bigl(B(x)+D(x)\bigr)\mathcal{P}(x;t)+B(x-1)\mathcal{P}(x-1;t)+D(x+1)\mathcal{P}(x+1;t),
\label{BDeq}
\end{align}
with the boundary condition(s)
\begin{equation}
D(0)=0,\quad B(N)=0:\  (\text{only  for a finite case}).
\label{bdcond1}
\end{equation}
The $(N+1)\times(N+1)$ matrix $L_{BD}$ is  tri-diagonal,
\begin{align}
&{L_{BD}}_{x+1\,x}=B(x),\ \  {L_{BD}}_{x-1\,x}=D(x),\ {L_{BD}}_{x\,x}=-B(x)-D(x), \n
&\hspace{30mm}{L_{BD}}_{x\,y}=0,\quad |x-y|\ge2,
\label{LBDdef}
\end{align}
representing the nearest neighbour interactions.

There are many exactly solvable cases with the explicit forms of the complete sets of eigenvalues and eigenfunctions of $L_{BD}$.
 Among them, those related with the hypergeometric orthogonal polynomials of Askey scheme
\cite{askey, ismail, koeswart, gasper} are pertinent for the present purpose.
All the discrete orthogonal polynomials of Askey scheme $\{\check{P}_n(x)\}$ with the normalised orthogonality measure
function $W(x)>0$,
\begin{equation}
\sum_{x\in\mathcal{X}}W(x)=1,\qquad \sum_{x\in\mathcal{X}}W(x)\check{P}_n(x)\check{P}_m(x)=0,\quad n\neq m,
\label{1orth}
\end{equation}
 satisfy a difference equation
 \begin{equation}
B(x)\!\left(\check{P}_n(x)-\check{P}_n(x+1)\right)\!+\!
D(x)\!\left(\check{P}_n(x)-\check{P}_n(x-1)\right)=\mathcal{E}(n)\check{P}_n(x),\ \mathcal{E}(n)\ge0,\  
n\in\mathcal{X},
\label{bdeq0}
\end{equation}
with positive coefficient functions $B(x)$ and $D(x)$ satisfying the boundary conditions \eqref{bdcond1}.
The orthogonality measure function $W(x)$ is determined by $B(x)$ and $D(x)$,
\begin{equation}
W(x)=W(0)\prod_{y=0}^{x-1}\frac{B(y)}{D(y+1)},\quad 1\le x\in\mathcal{X},
\label{BDW}
\end{equation}
in which the constant $W(0)$ is fixed by the normalisation condition $\sum_{x\in\mathcal{X}}W(x)=1$.
The following Theorem was shown in \cite{bdsol}.
%%%%%%%%%%%%%%%%%%%%%%%%%%%%%%%%%%%%%%%%%%%%%%%%%%%%%%
\begin{theo}
\label{theo1}
If the coefficient functions $B(x)$ and $D(x)$ of the Askey scheme polynomials \eqref{bdeq0} are 
chosen as the birth and death rates of $L_{BD}$ \eqref{LBDdef}, the corresponding
$\{\check{P}_n(x)\}$ and $W(x)$ provide the complete set of eigensystem of $L_{BD}$,
\begin{equation}
L_{BD}\check{P}_n(x)W(x)=-\mathcal{E}(n)\check{P}_n(x)W(x),\quad n\in\mathcal{X},
\label{LBDsol1}
\end{equation}
with $W(x)$ providing the stationary distribution with $\check{P}_0(x)\equiv1$,
\begin{equation}
\lim_{t\to+\infty}\mathcal{P}(x;t)=W(x), \quad \sum_{x\in\mathcal{X}}W(x)=1,\quad L_{BD}W(x)=0,
\label{stationary}
\end{equation}
as the zero-mode of $L_{BD}$,  i.e. $\mathcal{E}(0)=0$.
\end{theo}
%%%%%%%%%%%%%%%%%%%%%%%%%%%%%%%%%%%%%%%%%%%%%%%%%%%%
\begin{rema}
\label{askeypoly}
All of the discrete orthogonal polynomials of Askey scheme, ($q$-)Krawtchouk, (dual, $q$-)Hahn,
($q$-)Racah, etc could be named as Birth and Death polynomials. Those having the Jackson integral measures,
 e.g. the big $q$-Jacobi polynomials etc require slightly different formulations as detailed in {\rm \cite{os34}}.
 For the concrete expressions of the birth and death rates $B(x)$ and $D(x)$ for each polynomial, consult
{\rm \cite{bdsol, os12}}.
\end{rema}
%%%%%%%%%%%%%%%%%%%%%%%%%%%%%%%%%%%%%%%%%%%%%%%%%%%%%
%                                                   %
% 2.2 General setting                               %
%                                                   %
%%%%%%%%%%%%%%%%%%%%%%%%%%%%%%%%%%%%%%%%%%%%%%%%%%%%%
\subsection{Multivariable cases: General setting}
\label{sec:genset}

Let us now try and formulate the birth and death processes for $n$ groups of population
denoted by $\bm{x}=(x_1,x_2,\ldots,x_n)\in\mathcal{X}$.
The population space $\mathcal{X}$ is now a subset of ${\mathbb N}_0^n$.
The natural generalisation of the birth and death (BD) equation of  
the single species case \eqref{bdeqformal},\eqref{BDeq}, 
%the time evolution of the
%equation for probability distribution $\mathcal{P}(\bm{x};t)$ takes the following form 
reads as follows
\begin{align}
\frac{\partial}{\partial t}\mathcal{P}(\bm{x};t)&=(L_{BD}\mathcal{P})(\bm{x};t)
=\sum_{\bm{y}\in\mathcal{X}}{L_{BD}}_{\bm{x}\,\bm{y}}\mathcal{P}(\bm{y};t),
\quad \mathcal{P}(\bm{x};t)\ge0,\quad \sum_{\bm{x}\in\mathcal{X}} \mathcal{P}(\bm{x};t)=1,
\label{bdeqformal2}\\
&=-\sum_{j=1}^n(B_j(\bm{x})+D_j(\bm{x}))\mathcal{P}(\bm{x};t)
+\sum_{j=1}^nB_j(\bm{x}-\bm{e}_j)\mathcal{P}(\bm{x}-\bm{e}_j;t)\n
&\hspace{5.8cm}+\sum_{j=1}^nD_j(\bm{x}+\bm{e}_j)\mathcal{P}(\bm{x}+\bm{e}_j;t),
\label{BDeq2}
\end{align}
in which $\bm{e}_j$ is the  $j$-th unit vector, $j=1,\ldots,n$. 
The birth and death rates for $n$ groups $B_j(\bm{x})$, $D_j(\bm{x})$  are all positive with the
boundary conditions
\begin{equation}
\left.
\begin{array}{ccc}
B_j(\bm{x})=0  &  \text{if}\ \ \bm{x}\in\mathcal{X} & \text{and}\ \ \bm{x}+\bm{e}_j\notin\mathcal{X}   \\[2pt]
D_j(\bm{x})=0  &   \text{if}\ \ \bm{x}\in\mathcal{X}& \text{and}\ \ \bm{x}-\bm{e}_j\notin\mathcal{X}  
\end{array}
\right\}, \quad j=1,\ldots,n.
\label{becond2}
\end{equation}
This is a typical example of the nearest neighbour interactions in $n$ dimensions.

Let us impose the condition that the $|\mathcal{X}|\times|\mathcal{X}|$ matrix $L_{BD}$ has a stationary distribution
\begin{equation}
(L_{BD}W)(\bm{x})=0,\quad \sum_{\bm{x}\in\mathcal{X}}W(\bm{x})=1, \quad W(\bm{x})>0,
\quad \bm{x}\in\mathcal{X},
\label{stationaryn}
\end{equation}
which constrain $\{B_j(\bm{x}),D_j(\bm{x})\}$ severely.
In order to relate $W(\bm{x})$ with the birth/death rates, let us introduce the following notational conventions.
For a smooth function $f(x)$, the Taylor expansion reads
\begin{align*}
f(x+a)=\sum_{n=0}^\infty \frac{a^n}{n!}\frac{d^n}{dx^n}f(x)=(e^{a\partial}f)(x),\quad \partial\eqdef \frac{d}{dx}
\quad \Longrightarrow f(x\pm1)=(e^{\pm\partial}f)(x).
\end{align*}
Since only polynomials and the corresponding weight functions appear in this paper, it is quite safe to introduce the
operators $e^{\pm\partial_j}$ acting on a function $f(\bm{x})$ over $\mathcal{X}$,
\begin{equation*}
e^{\pm\partial_j}f(\bm{x})=f(\bm{x}\pm\bm{e}_j)\,e^{\pm\partial_j},\quad \partial_j\eqdef\frac{\partial}{\partial x_j},
\quad j=1,\ldots,n.
\end{equation*}
They are upper and lower triangular matrices acting on  vectors (functions) over $\mathcal{X}$
and transpose of each other
\begin{equation*}
\bigl(e^{\pm\partial_j}\bigr)_{\bm{x}\,\bm{y}}=\delta_{\bm{x}\pm\bm{e}_j\,\bm{y}},
\qquad \bigl(e^{\partial_j}\bigr)^T=e^{-\partial_j},\quad j=1,\ldots,n.
\end{equation*}
It is well known that the basic relations of these operators, 
$e^{\partial_j}e^{-\partial_j}=I_d=e^{-\partial_j}e^{\partial_j}$ are
broken at certain boundaries  of $\mathcal{X}$. But the violation of the rules does not cause any problem 
when these operators are used in conjunction with $\{B_j(\bm{x}),D_j(\bm{x})\}$ 
due to the boundary conditions \eqref{becond2}.
In terms of the operators $\{e^{\pm\partial_j}\}$,  
the birth and death operator $L_{BD}$ can be expressed succinctly as
\begin{align}
L_{BD}&=-\sum_{j=1}^n\left[B_j(\bm{x})-B_j(\bm{x}-\bm{e}_j)e^{-\partial_j}+D_j(\bm{x})-D_j(\bm{x}+\bm{e}_j)e^{\partial_j}\right],
\label{LBDop0}\\
&=-\sum_{j=1}^n\left[(1-e^{-\partial_j})B_j(\bm{x})+(1-e^{\partial_j})D_j(\bm{x})\right],
\n
&=-\sum_{j=1}^n(1-e^{-\partial_j})\bigl(B_j(\bm{x})-D_j(\bm{x}+\bm{e}_j)\,e^{\partial_j}\bigr).
\nonumber
\end{align}
This provides a sufficient condition for the zero mode of $L_{BD}$ \eqref{stationaryn}, the stationary distribution $W(\bm{x})>0$,
\begin{align}
\bigl(B_j(\bm{x})-D_j(\bm{x}+\bm{e}_j)e^{\partial_j}\bigr)W(\bm{x})=0
\ \Rightarrow \frac{W(\bm{x}+\bm{e}_j)}{W(\bm{x})}=\frac{B_j(\bm{x})}{D_j(\bm{x}+\bm{e}_j)},
\quad j=1,\ldots,n.
\label{zerocond}
\end{align}
These conditions determine entire $W(\bm{x})$ starting from the origin $W(\bm{0})$ so long as the 
two different routes give the same results,
\begin{align*}
W(\bm{x})\ &\to W(\bm{x}+\bm{e}_j)\ \to W(\bm{x}+\bm{e}_j+\bm{e}_k)\\[-4pt]
& \hspace{4.8cm} ||\\[-4pt]
& \to W(\bm{x}+\bm{e}_k)\ \to W(\bm{x}+\bm{e}_k+\bm{e}_j).
\end{align*}
%%%%%%%%%%%%%%%%%%%%%%%%%%%%%%%%%%
\begin{prop}
\label{compcond}
The existence of the stationary distribution imposes severe conditions on the  birth/death rates.
The sufficient condition that the BD process $L_{BD}$ \eqref{BDeq2} defined on a finite $\mathcal{X}$ 
to have a stationary distribution $W(\bm{x})$ is the following compatibility condition,
\begin{equation}
\frac{B_j(\bm{x})}{D_j(\bm{x}+\bm{e}_j)}\frac{B_k(\bm{x}+\bm{e}_j)}{D_k(\bm{x}+\bm{e}_j+\bm{e}_k)}
=\frac{B_k(\bm{x})}{D_k(\bm{x}+\bm{e}_k)}\frac{B_j(\bm{x}+\bm{e}_k)}{D_j(\bm{x}+\bm{e}_k+\bm{e}_j)},
\quad j,k=1,\ldots,n.
\end{equation}
For an infinite $\mathcal{X}$, the summability $\sum_{x\in\mathcal{X}}W(\bm{x})<\infty$ must be satisfied.
\end{prop}
%%%%%%%%%%%%%%%%%%%%%%%%%%%%%%%%%%%
Compatibility conditions for difference operators acting on various multivariable orthogonal polynomials 
are discussed in \cite{illxu}.
%%%%%%%%%%%%%%%%%%%%%%%%%%%%%%%%%%%%%%%%
\begin{rema}
\label{moreinfo}
The orthogonality weight function (stationary distribution) is determined by the birth/death rates.
That is, the birth/death rates have more information than the orthogonality measure.
\end{rema}
%%%%%%%%%%%%%%%%%%%%%%%%%%%%%%%%%%%%%%%%%%%
Supposing the stationary distribution  $W(\bm{x})$ is obtained, let us introduce a new operator $\mathcal{H}$
in terms of a similarity transformation of $L_{BD}$ by the square root of $W(\bm{x})$
\begin{equation}
\mathcal{H}\eqdef -\bigl(\sqrt{W(\bm{x})}\bigr)^{-1}L_{BD}\sqrt{W(\bm{x})}.
\label{Hdef}
\end{equation} \goodbreak
From  \eqref{LBDop0}, it is easy to see
\begin{align*}
&\bigl(\sqrt{W(\bm{x})}\bigr)^{-1}\bigl(B_j(\bm{x})+D_j(\bm{x})\bigr)\sqrt{W(\bm{x})}=B_j(\bm{x})+D_j(\bm{x}),\\
&\bigl(\sqrt{W(\bm{x})}\bigr)^{-1}D_j(\bm{x}+\bm{e}_j)\,e^{\partial_j}\sqrt{W(\bm{x})}\\
&\qquad =\bigl(\sqrt{W(\bm{x})}\bigr)^{-1}\sqrt{D_j(\bm{x}+\bm{e}_j)}\sqrt{D_j(\bm{x}+\bm{e}_j)}\sqrt{W(\bm{x}+\bm{e}_j)}\,e^{\partial_j}\\
& \qquad =\bigl(\sqrt{W(\bm{x})}\bigr)^{-1}\sqrt{D_j(\bm{x}+\bm{e}_j)}\sqrt{W(\bm{x})}\sqrt{B_j(\bm{x})}\,e^{\partial_j}
=\sqrt{B_j(\bm{x})D_j(\bm{x}+\bm{e}_j)}\,e^{\partial_j},\\
&\bigl(\sqrt{W(\bm{x})}\bigr)^{-1}B_j(\bm{x}-\bm{e}_j)\,e^{-\partial_j}\sqrt{W(\bm{x})}\\
&
\qquad =\bigl(\sqrt{W(\bm{x})}\bigr)^{-1}\sqrt{B_j(\bm{x}-\bm{e}_j)}\sqrt{B_j(\bm{x}-\bm{e}_j)}\sqrt{W(\bm{x}-\bm{e}_j)}\,e^{-\partial_j}\\
&\qquad =\bigl(\sqrt{W(\bm{x})}\bigr)^{-1}\sqrt{B_j(\bm{x}-\bm{e}_j)}\sqrt{W(\bm{x})}\sqrt{D_j(\bm{x})}\,e^{-\partial_j}
 =\sqrt{B_j(\bm{x}-\bm{e}_j)D_j(\bm{x})}\,e^{-\partial_j}.
\end{align*}
These lead to an expression of $\mathcal{H}$ depending only on the birth and death rates 
$\{B_j(\bm{x}),D_j(\bm{x})\}$,
\begin{align}
\mathcal{H}&=\sum_{j=1}^n\left[B_j(\bm{x})+D_j(\bm{x})-\sqrt{B_j(\bm{x})D_j(\bm{x}+\bm{e}_j)}\,e^{\partial_j}
-\sqrt{B_j(\bm{x}-\bm{e}_j)D_j(\bm{x})}\,e^{-\partial_j}\right],
\label{Hdef2}\\
\mathcal{H}_{\bm{x}\,\bm{y}}&=\sum_{j=1}^n\left[\bigl(B_j(\bm{x})+D_j(\bm{x})\bigr)\,\delta_{\bm{x}\,\bm{y}}-\sqrt{B_j(\bm{x})D_j(\bm{x}+\bm{e}_j)}\,\delta_{\bm{x}+\bm{e}_j\,\bm{y}}\right.\n
&\left.\hspace{5.2cm}
-\sqrt{B_j(\bm{x}-\bm{e}_j)D_j(\bm{x})}\,\delta_{\bm{x}-\bm{e}_j\,\bm{y}}\right],
\label{Hdef3}
\end{align}
ushering the following
%%%%%%%%%%%%%%%%%%%%%%%%%%%%%%%%%%%%%%%%%%%%%%%%%%%
\begin{prop}
\label{semidef}
The operator $\mathcal{H}$ is a real symmetric and positive semi-definite matrix over $\mathcal{X}$ as evidenced by
the factorisation form
\begin{align}
&\hspace{4cm}\mathcal{H}=\sum_{j=1}^n\mathcal{A}_j(\bm{x})^T\mathcal{A}_j(\bm{x}),
\qquad \mathcal{H}_{\bm{x}\,\bm{y}}=\mathcal{H}_{\bm{y}\,\bm{x}},
\label{Hfac}\\
&\mathcal{A}_j(\bm{x})\eqdef \sqrt{B_j(\bm{x})}-e^{\partial_j}\sqrt{D_j(\bm{x})},\ 
\mathcal{A}_j(\bm{x})^T= \sqrt{B_j(\bm{x})}-\sqrt{D_j(\bm{x})}\,e^{-\partial_j},\ 
j=1,\ldots,n.
\label{Ajdef}
\end{align}
Obviously $\sqrt{W(\bm{x})}$ is the zero mode of $\mathcal{A}_j(\bm{x})$ and $\mathcal{H}$
\begin{equation}
\mathcal{A}_j(\bm{x})\sqrt{W(\bm{x})}=0, \quad j=1,\ldots,n \quad \Longrightarrow \mathcal{H}\sqrt{W(\bm{x})}=0,
\label{zeromodes}
\end{equation}
as \eqref{zerocond} shows
\begin{equation*}
0=\mathcal{A}_j(\bm{x})\sqrt{W(\bm{x})}=\sqrt{B_j(\bm{x})}\sqrt{W(\bm{x})}-\sqrt{D_j(\bm{x}+\bm{e}_j)}\sqrt{W(\bm{x}+\bm{e}_j)}.
\end{equation*}
The birth/death rates $\{B_j(\bm{x}),D_j(\bm{x})\}$ contain many parameters.  For the generic values of the
parameters the degeneracy of eigenvalues do not occur and the orthogonality of the eigenvectors is guaranteed.
The birth and death operator $L_{BD}$ has negative semi-definite eigenvalues, as the eigensystem of $L_{BD}$
is related to that of $\mathcal{H}$ by the similarity transformation in terms of $\sqrt{W(\bm{x})}$ \eqref{Hdef}.
\end{prop}
%%%%%%%%%%%%%%%%%%%%%%%%%%%%%%%%%%%%%%%%%%%%%%%%%

Let us introduce another operator $\widetilde{\mathcal H}$ by a similarity transformation of $\mathcal{H}$
in terms of the square root of the stationary distribution $W(\bm{x})$,
\begin{align}
\widetilde{\mathcal H}&\eqdef \bigl(\sqrt{W(\bm{x})}\bigr)^{-1}\mathcal{H}\sqrt{W(\bm{x})}.
\label{Hthdef}
\end{align}
It is interesting to note that
\begin{align*}
&\bigl(\sqrt{W(\bm{x})}\bigr)^{-1}\mathcal{A}_j(\bm{x})\sqrt{W(\bm{x})}\\
&\qquad =
\bigl(\sqrt{W(\bm{x})}\bigr)^{-1}\left(\sqrt{B_j(\bm{x})}\sqrt{W(\bm{x})}-
\sqrt{D_j(\bm{x}+\bm{e}_j)}\sqrt{W(\bm{x}+\bm{e}_j)}e^{\partial_j}\right)\\
&\qquad =\sqrt{B_j(\bm{x})}\Bigl(1-e^{\partial_j}\Bigr),\\
&\bigl(\sqrt{W(\bm{x})}\bigr)^{-1}\mathcal{A}_j(\bm{x})^T\sqrt{W(\bm{x})}\\
&\qquad 
=\bigl(\sqrt{W(\bm{x})}\bigr)^{-1}\left(\sqrt{B_j(\bm{x})}\sqrt{W(\bm{x})}-
\sqrt{D_j(\bm{x})}\sqrt{W(\bm{x}-\bm{e}_j)}e^{-\partial_j}\right)\\
&\qquad
=\left(\sqrt{B_j(\bm{x})}-\frac{D_j(\bm{x})}{\sqrt{B_j(\bm{x}-\bm{e}_j)}}e^{-\partial_j}\right),\\
&\bigl(\sqrt{W(\bm{x})}\bigr)^{-1}\mathcal{A}_j(\bm{x})^T\mathcal{A}_j(\bm{x})\sqrt{W(\bm{x})}\\
&\qquad=B_j(\bm{x})\bigl(1-e^{\partial_j}\bigr)+D_j(\bm{x})\bigl(1-e^{-\partial_j}\bigr).
\end{align*}
%%%%%%%%%%%%%%%%%%%%%%%%%%%%%%%%%%%%%%%%%%%%%%%%%%%%%
\begin{prop}
\label{nAskey}
{\bf Multivariable  difference equation of  Askey scheme polynomials}\\
The operator $\widetilde{\mathcal H}$
\begin{equation}
\widetilde{\mathcal H}=\sum_{j=1}^n\left[B_j(\bm{x})\bigl(1-e^{\partial_j}\bigr)+D_j(\bm{x})\bigl(1-e^{-\partial_j}\bigr)\right],
\label{Hthdef2}
\end{equation}
provides the left hand side of the multivariable generalisation of the difference equation \eqref{bdeq0} governing the Askey polynomials.
Obviously a constant is the zero mode of  $\widetilde{\mathcal H}$
\begin{equation}
\widetilde{\mathcal H}\,1=0,
\label{Hthzero}
\end{equation}
as is always the case with the polynomials in Askey scheme. For finite $|\mathcal{X}|$, the eigensystem of the 
real symmetric matrix $\mathcal{H}$ always exists and correspondingly the eigenvectors of $\widetilde{\mathcal H}$
also exist. But it is not evident in general that the eigenvectors can be expressed as polynomials in  
$\bm{x}$ or some simple functions of it.
\end{prop}
%%%%%%%%%%%%%%%%%%%%%%%%%%%%%%%%%%%%%%%%%%%%%%%%%%%%%%%
\begin{rema}
\label{rati}
In contrast to $\mathcal{H}$ \eqref{Hdef2},\eqref{Hdef3}, $\widetilde{\mathcal H}$ \eqref{Hthdef2} is square root free.
It is amenable for the determination of the explicit forms of the polynomials, as  shown in the next section.
In other words, the real symmetric and positive semi-definite matrix (linear operator) $\mathcal{H}$ which
governs the orthogonality, seems to have been put outside of the stage due to the square root factor.
\end{rema}
%%%%%%%%%%%%%%%%%%%%%%%%%%%%%
This program works well for a particular choice of the birth/death rates which  give rise to multivariate 
Krawtchouk polynomials as explored in the next section.

%%%%%%%%%%%%%%%%%%%%%%%%%%%%%%%%%%%%%
%
%         section3 Multivariate Krawtchouk polynomials %
%
%%%%%%%%%%%%%%%%%%%%%%%%%%%%%%%%%%%%%
\section{Multivariate Krawtchouk polynomials}
\label{sec:nKr}
Following the procedures introduced in the previous section, $n$-parameter dependent $n$-variate 
Krawtchouk polynomials are constructed explicitly based on a special choice of the birth/death rates.
For comparison, the basic data of the single variable Krawtchouk polynomials are recapitulated here \cite{koeswart}.
The birth and death rates are
\begin{align}
&B(x)=p(N-x),\quad  D(x)=(1-p)x,\quad  0<p<1,
\label{KRbdr}\\
& \Rightarrow
\widetilde{\mathcal H}=p(N-x)(1-e^\partial)+(1-p)(1-e^{-\partial}).
\label{KrHth}
\end{align}
The stationary distribution is the binomial distribution with the probability $p$,
\begin{equation}
W(p;x)=\binom{N}{x}p^x(1-p)^{N-x},\quad \mathcal{X}=\{0,1,\ldots,N\}.
\end{equation}
The eigenpolynomial of the difference equation is a truncated hypergeometric function ${}_2F_1$ with a
linear spectrum,
\begin{equation}
\widetilde{\mathcal H}P_n(x)=n\,P_n(x),
\qquad P_n(x)={}_2F_1\Bigl(\genfrac{}{}{0pt}{}{-n,\,-x}{-N}\Bigm|p^{-1}\Bigr)
=\sum_{k=0}^n\frac{(-n)_k(-x)_k}{(-N)_k}\frac{p^{-k}}{k!},
\label{Krform}
\end{equation}
in which $(a)_n$ is the shifted factorial defined for $a\in\mathbb{C}$ and nonnegative integer $n$,
$(a)_0=1$, $(a)_n=\prod_{k=0}^{n-1}(a+k)$, $n\ge1$.
The orthogonality reads
\begin{equation}
\sum_{x\in\mathcal{X}}W(p;x)P_n(x)P_m(x)=\frac{\delta_{n\,m}}{\binom{N}{n}\bigl(\frac{p}{1-p}\bigr)^n}.
\label{Krorth}
\end{equation}
The generating function of the Krawtchouk polynomials $\{P_n(x)\}$ is, see \cite{koeswart} (9.11.11),
\begin{equation}
G(p,x;t)\eqdef\left(1-\frac{(1-p)}{p}t\right)^x\bigl(1+t\bigr)^{N-x}=\sum_{n\in\mathcal{X}}\binom{N}{n}P_n(x)t^n.
\label{Krgen}
\end{equation}

%%%%%%%%%%%%%%%%%%%%%%%%%%%%%%%%%%%%%
%
%         section3.1  Main Results %
%
%%%%%%%%%%%%%%%%%%%%%%%%%%%%%%%%%%%%%
\subsection{Main Results}
\label{sec:Main}
Let us start with the summary as the following
\begin{theo}
\label{theo:main}
{\bf Main Theorem}\\
The birth and death rates are defined by two positive integers $N$ and $n$ ($N>n\ge2$) and $n$  positive numbers
$p_i>0$,  $i=1,\ldots, n$, 
\begin{equation}
B_i(\bm{x})\eqdef (N-|x|),\  \bm{x}=(x_1,\ldots,x_n)\in\mathbb{N}_0^n,\ |x|\eqdef \sum_{i=1}^nx_i,\
D_i(\bm{x})\eqdef p_i^{-1}x_i,\ i=1,\ldots,n.
\label{BDKrdef}
\end{equation}
The  population space is determined by the positivity of $\{B_i(\bm{x})\}$,
\begin{equation}
\mathcal{X}=\{\bm{x}\in\mathbb{N}_0^n\ |\, |x|\le N\},\quad |\mathcal{X}|=\binom{N+n}{n}.
\end{equation}
The compatibility condition \eqref{compcond} is trivially satisfied and  the stationary distribution is the multinomial distribution with probabilities $\{\eta_i\}$ which are functions of 
$\{p_i\}$,
\begin{align}
&W(\eta;\bm{x})=\frac{N!}{x_1!\cdots x_n!x_0!}\prod_{i=0}^n\eta_i^{x_i}
=\binom{N}{\bm{x}}\eta_0^{x_0}\bm{\eta}^{\bm{x}},\  x_0\eqdef
N-|x|,\  \binom{N}{\bm{x}}\eqdef\frac{N!}{x_1!\cdots x_n!x_0!},
\label{WKr}\\
&\quad \eta_i\eqdef\frac{p_i}{1+\sum_{j=1}^np_j},\quad \eta_0\eqdef\frac1{1+\sum_{i=1}^np_i},
\quad \sum_{i=0}^n\eta_i=1, \quad \bm{\eta}^{\bm{x}}\eqdef\prod_{i=1}^n\eta_i^{x_i}.
\label{etadef}
\end{align}
The eigenpolynomials $\{P_{\bm{m}}(\bm{x})\}$ of the difference operator $\widetilde{\mathcal H}$ \eqref{Hthdef2},
\begin{align}
&\widetilde{\mathcal H}=(N-|x|)\sum_{i=1}^n(1-e^{\partial_j})+\sum_{i=1}^np_i^{-1}x_i(1-e^{-\partial_j}),
\label{HthnKr}\\
&\widetilde{\mathcal H}P_{\bm{m}}(\bm{x})=\mathcal{E}(\bm{m})P_{\bm{m}}(\bm{x}),\qquad \bm{x},\bm{m}\in\mathcal{X}.
\label{Pmeigpol}
\end{align}
are indexed by the degrees $\bm{m}=(m_1,m_2,\ldots,m_n)\in\mathbb{N}_0^n$ in the dual  space
of $\mathcal{X}$, which is isomorphic to $\mathcal{X}$.
The eigenvalue $\mathcal{E}(\bm{m})$ has a linear spectrum
\begin{equation}
\mathcal{E}(\bm{m})\eqdef\sum_{j=1}^nm_j\lambda_j,
\label{eigform}
\end{equation}
in which $\lambda_j$ is the $j$-th root of a degree $n$ characteristic polynomial $\mathcal{F}(\lambda)$
of an $n\times n$ positive matrix $F(p)$  depending on $\{p_i\}$,
\begin{equation}
0=\mathcal{F}(\lambda)\eqdef Det\bigl(\lambda I_n-F(p)\bigr),\quad F(p)_{i\,j}\eqdef 1+p_i^{-1}\delta_{i\,j}.
\label{chareq}
\end{equation}
They are all positive, by construction.
The eigenpolynomial $P_{\bm{m}}(\bm{x})$ is a truncated   $(n+1,2n+2)$ hypergeometric function 
of Aomoto-Gelfand  {\rm \cite{AK,gelfand,mizu}}
\begin{gather}
\label{Pm}
P_{\bm{m}}(\bm{x})
\eqdef \sum_{\substack{\sum_{i,j}c_{ij}\leq N\\
(c_{ij})\in M_{n}({\mathbb N_{0}})}}
\frac{\prod\limits_{i=1}^{n}(-x_{i})_{\sum\limits_{j=1}^{n}c_{ij}}
\prod\limits_{j=1}^{n}(-m_{j})_{\sum\limits_{i=1}^{n}c_{ij}}}
{(-N)_{\sum_{i,j}c_{ij}}} \; \frac{\prod(u_{ij})^{c_{ij}}}{\prod c_{ij}!},
\end{gather}
in which $M_{n}({\mathbb N}_{0})$ is the set of square matrices of degree $n$ with nonnegative integer
elements. 
This looks like a simple generalisation of the single variable case \eqref{Krform} consisting of the
shifted factorials of the negative coordinates, degrees and N times the power of $u_{i\,j}$ divided by
the corresponding factorial.
The $n\times n$ matrix $u_{i\,j}$, which could be considered as the ``inverse probability" 
counterpart corresponding to the pair of 
the coordinate $x_i$ and the degree $m_j$,  carries the information of the system. It is defined by
\begin{equation}
u_{i\,j}\eqdef\frac{\lambda_j}{\lambda_j-p_i^{-1}},\quad i,j=1,\ldots,n.
\label{udef}
\end{equation}
The orthogonality of the eigenpolynomials of $\widetilde{\mathcal H}$ 
 is the  consequence  of the self-adjointness of $\mathcal{H}$ \eqref{Hfac} 
and the non-degenerary of eigenvalues for generic parameters. The norms of the eigenpolynomials 
are parametrised by the probability parameters $\{\bar{p}_j\}$ 
\begin{align}
\sum_{\bm{x}\in\mathcal{X}}W(\eta;\bm{x})P_{\bm{m}}(\bm{x})P_{\bm{m}'}(\bm{x})
&=\frac{\delta_{\bm{m}\,\bm{m}'}}{\binom{N}{\bm{m}}(\bar{\bm{p}})^{\bm{m}}},\qquad 
(\bar{\bm{p}})^{\bm{m}}\eqdef\prod_{j=1}^n\bar{p}_j^{m_j},
\label{nKrorth}\\
\bar{p}_j&=\Bigl(\sum_{i=1}^n\eta_iu_{i,j}^2-1\Bigr)^{-1}>0,\quad j=1,\ldots,n,
\label{pddef}
\end{align}
which  are positive, by construction.
The generating function of the above hypergeometric function \eqref{Pm} is well-known  {\rm\cite{mizu}},
%\begin{gather}\label{gen}
%\prod_{i=0}^{n}\left(\sum_{j=0}^{n}a_{ij}t_{j}\right)^{x_{i}}=\sum_{\bm{m} \in\mathcal{X}}
%\binom{N}{\bm{m}}P_{\bm{m}}(\bm{x})t_0^{m_0}\bm{t}^{\bm{m}},\quad
%\bm{t}^{\bm{m}}\eqdef\prod_{j=1}^nt_j^{m_j},\quad m_0\eqdef N-|m|,
%\end{gather}
\begin{align}
G(u,\bm{x};t)&\eqdef \prod_{i=0}^{n}\left(\sum_{j=0}^{n}b_{ij}t_{j}\right)^{x_{i}}=\sum_{\bm{m} \in\mathcal{X}}
\binom{N}{\bm{m}}P_{\bm{m}}(\bm{x})t_0^{m_0}\bm{t}^{\bm{m}},
\label{gen}\\
&\bm{t}^{\bm{m}}\eqdef\prod_{j=1}^nt_j^{m_j},\quad m_0\eqdef N-|m|.
\end{align}
The parameters $\{b_{i\,j}\}$ are related to $\{u_{i\,j}\}$ 
\begin{align}
b_{0\,j}=b_{i\,0}=1 \ {\rm for}  \ 0\leq i,j\leq n \  {\rm and} \ b_{i\,j}=1-u_{i\,j}\  {\rm for} \ i,j=1,\ldots,n.
\label{aurel}
\end{align}
%in which $a_{0\,j}=a_{i\,0}=1$ for $0\leq i,j\leq n$ and $a_{i\,j}=1-u_{i\,j}$ for $i,j=1,\ldots,n$.
\end{theo}

The rest of the section is devoted to the derivation of the main results, step by step.

%%%%%%%%%%%%%%%%%%%%%%%%%%%%%%%%%%%%%
%
%         section3.2  Stationary Distribution
%
%%%%%%%%%%%%%%%%%%%%%%%%%%%%%%%%%%%%%
\subsection{Stationary Distribution}
\label{sec:distri}
It is easy to verify that the compatibility condition \eqref{compcond} is trivially satisfied, as
\begin{equation*}
\frac{(N-|x|)p_j}{(x_j+1)}\frac{(N-|x|-1)p_k}{(x_k+1)}=
\frac{(N-|x|)p_k}{(x_k+1)}\frac{(N-|x|-1)p_j}{(x_j+1)}.
\end{equation*}
The stationary distribution can be obtained by using the two term elation 
$\frac{W(\bm{x}+\bm{e}_j)}{W(\bm{x})}=\frac{B_j(\bm{x})}{D_j(\bm{x}+\bm{e}_j)}$ \eqref{zerocond} 
starting from $W(\bm{0})$ step by step, which is the generalisation of the single variable formula \eqref{BDW}.
But it is easier to show that the stationary distribution \eqref{WKr}
\begin{equation*}
W(\eta;\bm{x})=\binom{N}{\bm{x}}\prod_{i=1}^np_i^{x_i}\eta_0^N,
\end{equation*}
satisfies the two term relation
\begin{align*}
\left((N-|x|)-p_j^{-1}(x_j+1)e^{\partial_j}\right)W(\eta;\bm{x})=0,\quad j=1,\ldots,n,
\end{align*}
as
\begin{align*}
e^{\partial_j}W(\eta;\bm{x})=\frac{N-|x|}{x_j+1}p_jW(\eta;\bm{x})
\ \Rightarrow p_j^{-1}(x_j+1)e^{\partial_j}W(\eta;\bm{x})=(N-|x|)W(\eta;\bm{x}).
\end{align*}

%%%%%%%%%%%%%%%%%%%%%%%%%%%%%%%%%%%%%
%
%         section3.3  Degree one Eigenpolynomials
%
%%%%%%%%%%%%%%%%%%%%%%%%%%%%%%%%%%%%%
\subsection{Degree one Eigenpolynomials}
\label{sec:onedegi}
Now the stationary distribution is established, the  resulting real symmetric and positive semi-definite
matrix (linear operator) $\mathcal{H}$ \eqref{Hdef2},\eqref{Hdef3},  and the difference equation operator 
$\widetilde{\mathcal H}$  \eqref{Hthdef2}
\begin{equation*}
\widetilde{\mathcal H}=\bigl(N-|x|)\sum_{i=1}^n(1-e^{\partial_j})
+\sum_{i=1}^np_i^{-1}x_i(1-e^{-\partial_j}),
\tag{\ref{HthnKr}}
\end{equation*}
can be employed to determine the explicit forms of the corresponding orthogonal polynomials.
%%%%%%%%%%%%%%%%%%%%%%%%%%%%%%%%%%%%%%%%%%%%%
\begin{prop}
\label{Hthinv}
It is obvious that the space of polynomials in $\bm{x}$ of maximal total degree $M$
\begin{equation}
V_M(\bm{x})={\rm Span}\left\{\bm{x}^{\bm{m}}\mid 0\le |m| \le M\right\},
\quad \bm{x}^{\bm{m}}\eqdef\prod_{i=1}^nx_i^{m_i}
\end{equation}
is invariant under the action of $\widetilde{\mathcal H}$
\begin{equation}
\widetilde{\mathcal H}V_M(\bm{x})\subseteq V_M(\bm{x}).
\end{equation}
There are $\binom{M+n-1}{n-1}$ eigenpolynomials of $\widetilde{\mathcal H}$ in $\bm{x}$ of maximal degree $M$.
\end{prop}
%%%%%%%%%%%%%%%%%%%%%%%%%%%%%%%%%%%%%%%%%%%%
Let us determine $n$ degree 1 eigenpolynomials of $\widetilde{\mathcal H}$ \eqref{HthnKr} with 
unknown coefficients $\{a_i\}$ and unit constant part
as is always the case with the hypergeometric polynomials of discrete variables,
\begin{align}
&P_{|m|=1}(\bm{x})=1+\sum_{i=1}^na_ix_i,
\quad
\widetilde{\mathcal H}P_{|m|=1}(\bm{x})=\lambda P_{|m|=1}(\bm{x}),
\label{deg1form}\\
&\Rightarrow -(N-|x|)\sum_{i=1}^na_i+\sum_{i=1}^np_i^{-1}a_ix_i=\lambda\Bigl(1+\sum_{i=1}^na_ix_i\Bigr).
\nonumber
\end{align}
By equating the coefficients of $x_i$, an eigenvalue equations of $\{a_i\}$ are obtained,
\begin{align}
\sum_{k=1}^na_k+p_i^{-1}a_i&=\lambda a_i,\quad i=1,\ldots,n,
\label{1eig}\\
-N\sum_{k=1}^na_k&=\lambda.
\label{lambdaeq}
\end{align}
Then $\lambda$ is the  root of a degree $n$ characteristic polynomial $\mathcal{F}(\lambda)$
of an $n\times n$ positive matrix $F(p)$  depending on $\{p_i\}$,
\begin{equation*}
0=\mathcal{F}(\lambda)\eqdef \text{Det}\bigl(\lambda I_n-F(p)\bigr),\quad F(p)_{i\,j}\eqdef 1+p_i^{-1}\delta_{i\,j}.
\tag{\ref{chareq}}
\end{equation*}
For each eigenvalue $\lambda_j$, satisfying the relation
\begin{equation}
\sum_{i=1}^n\frac{1}{\lambda_j-p_i^{-1}}\equiv\sum_{i=1}^n\frac{p_i}{p_i\lambda_j-1}=1,\quad j=1,\ldots,n,
\label{plambrel}
\end{equation}
the unknown coefficients $\{a_i\}$ are determined,
\begin{equation*}
a_{i,j}=-\frac{\lambda_j}{N(\lambda_j-p_i^{-1})},\quad i,j=1,\ldots,n.
\end{equation*}
Let us tentatively identify the above $j$-th solution as $\bm{m}=\bm{e}_j$ solution
\begin{equation}
P_{\bm{e}_j}(\bm{x})=1-\frac1N\sum_{i=1}^n\frac{\lambda_j}{\lambda_j-p_i^{-1}}x_i,\quad j=1,\ldots,n.
\label{e_jsol}
\end{equation}
By comparing these polynomials with the general hypergeometric functions \cite{AK, gelfand}
and, in particular, that corresponding to  the generalised form of $n$-variate Krawtchouk polynomials \cite{mizu}
\begin{gather*}
P_{\bm{m}}(\bm{x})
\eqdef \sum_{\substack{\sum_{i,j}c_{ij}\leq N\\
(c_{ij})\in M_{n}({\mathbb N_{0}})}}
\frac{\prod\limits_{i=1}^{n}(-x_{i})_{\sum\limits_{j=1}^{n}c_{ij}}
\prod\limits_{j=1}^{n}(-m_{j})_{\sum\limits_{i=1}^{n}c_{ij}}}
{(-N)_{\sum_{i,j}c_{ij}}} \; \frac{\prod(u_{ij})^{c_{ij}}}{\prod c_{ij}!},
\tag{\ref{Pm}}
\end{gather*}
the system parameters $\{u_{i\,j}\}$ are completely identified
\begin{equation*}
P_{\bm{e}_j}(\bm{x})=1-\frac1N\sum_{i=1}^nu_{i\,j}x_i,\qquad 
u_{i\,j}=\frac{\lambda_j}{\lambda_j-p_i^{-1}},\quad i,j=1,\ldots,n.
\tag{\ref{udef}}
\end{equation*}
%%%%%%%%%%%%%%%%%%%%%%%%%%%%%%%%%%%%%%%%%%%%
\begin{prop}
\label{1all}
When the particular hypergeometric form \eqref{Pm} 
 is assumed for the eigenpolynomials of $\widetilde{\mathcal H}$ \eqref{HthnKr}, the explicit forms of the degree1
 polynomials determine the entire set of the polynomials.
\end{prop}
%%%%%%%%%%%%%%%%%%%%%%%%%%%%%%%%%%%%%%%%%%%%
\begin{rema}
\label{sinusoidal}
If we write the lowest degree  discrete orthogonal polynomial of Askey schme
\begin{equation*}
\check{P}_1(x)=1+const. \zeta(x),\quad \zeta(0)=0,
\end{equation*}
the higher polynomials are the polynomials in $\zeta(x)$. That is,  they  
are also obtained by Gram-Schmidt orthonormalisation with the orthogonality measure 
given by the birth/death rates. The scale of $\zeta(x)$ is immaterial. 
They are, for example, $\zeta(x)=x,\ x(x+d),\ q^{-x}-1,\ 1-q^x,\ (q^{-x}-1)(1-d q^x)$. 
These $\zeta(x)$ have many special properties and they are called `sinusoidal coordinates',
{\rm\cite{os12}}.
\end{rema}
%%%%%%%%%%%%%%%%%%%%%%%%%%%%%%%%%%%%%%%%%%%%
%%%%%%%%%%%%%%%%%%%%%%%%%%%%%%%%%%%%%
%
%         section3.4  General Eigenpolynomials
%
%%%%%%%%%%%%%%%%%%%%%%%%%%%%%%%%%%%%%
\subsection{General Eigenpolynomials}
\label{sec:genpol}
The next task is to verify that the higher degree ones $\{P_{\bm{m}}(\bm{x})\}$ \eqref{Pm} are the eigenpolynomials of 
$\widetilde{\mathcal H}$ \eqref{HthnKr}, too. For this,  the explicit forms of the eigenvalues $\mathcal{E}(\bm{m})$ 
\eqref{eigform} are necessary. Since $P_{\bm{m}}(\bm{x})$ has the form
\begin{equation}
P_{\bm{m}}(\bm{x})=1-\frac1N\sum_{i,j=1}^nx_im_ju_{i\,j}+ \text{higher degrees},
\end{equation}
 $\widetilde{\mathcal H}$ acting on the higher degrees produces only the terms of linear and higher degrees.
The only constant part of $\widetilde{\mathcal H}P_{\bm{m}}(\bm{x})$ comes from 
$N\sum_{i=1}^n(1-e^{\partial_j})$ acting on the linear part,
\begin{equation*}
N\sum_{k=1}^n(1-e^{\partial_k})\left\{-\frac1N\sum_{i,j=1}^nx_im_ju_{i\,j}\right\}
=\sum_{i\,j}m_ju_{i\,j}=\sum_{j=1}^nm_j\lambda_j\sum_{i=1}^n\frac{1}{\lambda_j-p_i^{-1}}
=\sum_{j=1}^nm_j\lambda_j,
\end{equation*}
in which \eqref{plambrel} is used.  
After applying $(1-e^{\partial_i})$, all higher degree terms vanish at the origin ${\bm x}={\bm 0}$ as they
consist of terms like $(x_i)_k(x_j)_l$, $k+l\ge2$.
This leads to the following
%%%%%%%%%%%%%%%%%%%%%%%%%%%%%%%%%%%%%%%%%%%%
\begin{prop}
\label{lineareig}
If $P_{\bm{m}}(\bm{x})$ \eqref{Pm} is an eigenpolynomial of $\widetilde{\mathcal H}$ \eqref{HthnKr},
it has a linear spectrum
\begin{equation*}
\mathcal{E}(\bm{m})\eqdef\sum_{j=1}^nm_j\lambda_j.
\tag{\ref{eigform}}
\end{equation*}
\end{prop}
%%%%%%%%%%%%%%%%%%%%%%%%%%%%%%%%%%%%%%%%%%%%
The next task is to prove
\begin{equation*}
\widetilde{\mathcal H}P_{\bm{m}}(\bm{x})
=\Bigl(\sum_{k=1}^nm_k\lambda_k\Bigr)P_{\bm{m}}(\bm{x}),\qquad \bm{x},\bm{m}\in\mathcal{X}.
\tag{\ref{Pmeigpol}}
\end{equation*}
It is based on the generating function \eqref{gen}.  As for the single variable Krawtchouk polynomial,
it is straightforward to show with $\widetilde{\mathcal H}$ \eqref{KrHth} and the generating function
$G(p,x;t)$ \eqref{Krgen} that
\begin{align*}
\widetilde{\mathcal H}G(p,x;t)=t\frac{\partial}{\partial t}G(p,x;t)&=\sum_{n=0}^N\binom{N}{n}nP_n(x)t^n,\\
\Rightarrow \widetilde{\mathcal H}P_n(x)&=nP_n(x).
\end{align*}
The corresponding formula for the $n$-variable generating function $G(u,\bm{x},t)$ \eqref{gen} reads
\begin{equation}
\widetilde{\mathcal H}G(u,\bm{x};t)
=\left(\sum_{k=1}^n\lambda_kt_k\frac{\partial}{\partial t_k}\right)G(u,\bm{x};t),
\label{genform}
\end{equation}
which leads to \eqref{Pmeigpol} above. A bit lengthy but straightforward derivation of \eqref{genform} is listed below,
as it is the only substantial calculation in this paper. 
The generating function  $G(u,\bm{x},t)$ \eqref{gen} is decomposed as
\begin{equation*}
G(u,\bm{x},t)=\bigl(t_0+|t|\bigr)^{x_0}\prod_{i=1}^nT_i^{x_i},\quad |t|\eqdef \sum_{k=1}^nt_k,
\quad T_i\eqdef t_0+\sum_{j=1}^nb_{i\,j}t_j.
\end{equation*}
The action of the two parts of $\widetilde{\mathcal H}$
\begin{equation*}
(N-|x|)\sum_{i=1}^n(1-e^{\partial_i})=x_0\sum_{i=1}^n(1-e^{\partial_i}),\qquad
\sum_{i=1}^np_i^{-1}x_i(1-e^{-\partial_i}),
\end{equation*}
on $G(u,\bm{x},t)$ is evaluated separately. The first part gives
\begin{align*}
x_0\sum_{i=1}^n(1-e^{\partial_i})G(u,\bm{x};t)&=
x_0\bigl(t_0+|t|\bigr)^{x_0-1}\prod_{l=1}^nT_l^{x_l}\times
\sum_{i=1}^n\sum_{k=1}^n(1-b_{i\,k})t_k\\
&=\sum_{k=1}^n\lambda_kt_k\,x_0\bigl(t_0+|t|\bigr)^{x_0-1}\prod_{l=1}^nT_l^{x_l},
\qquad \qquad  (*)
\end{align*}
in which \eqref{aurel}, \eqref{udef} and \eqref{plambrel} are used to obtain
\begin{equation*}
\sum_{i=1}^n\sum_{k=1}^n(1-b_{i\,k})t_k=\sum_{i=1}^n\sum_{k=1}^nu_{i\,k}t_k
=\sum_{i=1}^n\sum_{k=1}^n\frac{\lambda_k}{\lambda_k-p_i^{-1}}t_k=\sum_{k=1}^n\lambda_kt_k.
\end{equation*}
The second part gives
\begin{align*}
\sum_{i=1}^np_i^{-1}x_i(1-e^{-\partial_i})G(u,\bm{x};t)
&=\sum_{i=1}^np_i^{-1}x_i\bigl(t_0+|t|\bigr)^{x_0}\prod_{l\neq i}^nT_l^{x_l}
\cdot T_i^{x_i-1}\sum_{k=1}^n(-1)(1-b_{i\,k})t_k\\
&=\sum_{i=1}^n\sum_{k=1}^n\lambda_kb_{i\,k}t_k x_i\bigl(t_0+|t|\bigr)^{x_0}\prod_{l\neq i}^nT_l^{x_l}
\cdot T_i^{x_i-1}, \qquad (**)
\end{align*}
in which
\begin{equation*}
p_i^{-1}(-1)(1-b_{i\,k})t_k=-p_i^{-1}\frac{\lambda_k}{\lambda_k-p_i^{-1}}t_k=\lambda_kb_{i\,k}t_k,\quad 
b_{i\,k}=1-\frac{\lambda_k}{\lambda_k-p_i^{-1}}=-\frac{p_i^{-1}}{\lambda_k-p_i^{-1}},
\end{equation*}
are used. Now r.h.s. of \eqref{genform} reads
\begin{align}
&\left(\sum_{k=1}^n\lambda_kt_k\frac{\partial}{\partial t_k}\right)G(u,\bm{x};t)\n
&\qquad=x_0\sum_{k=1}^n\lambda_kt_k\ \bigl(t_0+|t|\bigr)^{x_0-1}\prod_{l=1}^nT_l^{x_l}\n
&\qquad \quad + \sum_{i=1}^n\sum_{k=1}^n\lambda_kb_{i\,k}t_k x_i\bigl(t_0+|t|\bigr)^{x_0}\prod_{l\neq i}^nT_l^{x_l}
\cdot T_i^{x_i-1}.
\end{align}
The r.h.s. are qual to $(*)+(**)$ and this concludes the proof, which leads to the following
%%%%%%%%%%%%%%%%%%%%%%%%%%%%%%%%%%%%%%%%%%%%%%%%%%%%
\begin{prop}
\label{eigproof}
The multivariate Krawtchouk polynomials  $\{P_{\bm{m}}(\bm{x})\}$\! \eqref{Pm} constitute the complete 
set of eigenpolynomils of the difference operator $\widetilde{\mathcal H}$ \eqref{HthnKr} which is derived from
the $n$-variate Birth and Death process \eqref{BDKrdef}.
\end{prop}
%%%%%%%%%%%%%%%%%%%%%%%%%%%%%%%%%%%%%%%%%%%%%%%%%%%%%
\begin{prop}{\bf $\mathfrak{S}_n$ Symmetry}
\label{Snsymm}
The multivariate Krawtchouk polynomials  $\{P_{\bm{m}}(\bm{x})\}$\! \eqref{Pm} is 
invariant under the symmetric group $\mathfrak{S}_n$, 
due to the arbitrariness of the ordering of $n$ roots $\{\lambda_j\}$ 
of the characteristic equation \eqref{chareq} in the parameters $u_{i\,j}$ \eqref{udef}.
\end{prop}

%%%%%%%%%%%%%%%%%%%%%%%%%%%%%%%%%%%%%
%
%         section3.5  Orthogonality
%
%%%%%%%%%%%%%%%%%%%%%%%%%%%%%%%%%%%%%
\subsection{Orthogonality revisited}
\label{sec:ortho}
The orthogonality of the generic truncated hypergeometric function \eqref{Pm}
has been discussed by many authors \cite{Gri1, gr1, mizu, diaconis13, Gri2}. 
Among them, an interpretation of Mizukawa's result \cite{mizu} is  presented from the view point of 
this paper. Now  a half of his result with my reinterpretation is recapitulated  as the following
%%%%%%%%%%%%%%%%%%%%%%%%%%%%%%%%%%%%%%%%%%%%%%%%%%
\begin{prop}
\label{mizuT}
{\bf [Mizukawa \cite{mizu} Theorem 1]}\\
If a diagonal matrix $D_1=diag(\eta_0,\eta_1,\ldots,\eta_n)$ is transformed to another 
diagonal matrix $D_2=diag(1,\bar{p}_1^{-1},\ldots,\bar{p}_n^{-1})$,
by an $(n+1)\times(n+1)$
matrix $B=\{b_{i\,j}\}$ defined in \eqref{aurel}, 
\begin{equation*}
B^TD_1B=D_2,
\end{equation*}
then the generic hypergeometric function 
$P_{\bm{m}}(\bm{x})$ \eqref{Pm} are orthogonal with respect to the multinomial distribution
$W(\eta,\bm{x})$ \eqref{WKr}
%\begin{equation*}
%\sum_{\bm{x}\in\mathcal{X}}W(\eta;\bm{x})P_{\bm{m}}(\bm{x})P_{\bm{m}'}(\bm{x})
%=\frac{\delta_{\bm{m}\,\bm{m}'}}{\binom{N}{\bm{m}}\bar{\bm{\eta}}^{\bm{m}}},\qquad 
%\bar{\bm{\eta}}^{\bm{m}}\eqdef\prod_{j=1}^n\bar{\eta}_j^{m_j}.
%\tag{\ref{nKrorth}}
%\end{equation*}
\begin{equation*}
\langle P_{\bm{m}},P_{\bm{m}'}\rangle=
\frac{\delta_{\bm{m}\,\bm{m}'}}{\binom{N}{\bm{m}}\bar{\bm{p}}^{\bm{m}}},\qquad 
\bar{\bm{p}}^{\bm{m}}\eqdef\prod_{j=1}^n\bar{p}_j^{m_j}.
\tag{\ref{nKrorth}}
\end{equation*}
Here the notation
\begin{equation*}
\langle f,g\rangle\eqdef \sum_{\bm{x}\in\mathcal{X}}W(\eta;\bm{x})f(\bm{x})g(\bm{x}),
\end{equation*}
is used.
It should be stressed that $\{u_{i\,j}\}$ and $\{\eta_i\}$ are completely generic in this context.
\end{prop}
%%%%%%%%%%%%%%%%%%%%%%%%%%%%%%%%%%%%%%%%%%%%%%%
It would be easier to understand and remember the above Proposition \ref{mizuT}, if it is rephrased by 
the following
%%%%%%%%%%%%%%%%%%%%%%%%%%%%%%%%%%%%%%%%%%%%%%%%%
\begin{rema}
\label{rephmizu}
The orthogonality \eqref{nKrorth} holds if the degree one polynomials 
$\{P_{\bm{e}_j}(\bm{x})\}$ are orthogonal to 1 and with each other,
\begin{align}
\langle 1,P_{\bm{e}_j}\rangle&=0\ \Longleftrightarrow \sum_{i=1}^n\eta_iu_{i\,j}=1,
\qquad j=1,\ldots,n,
\label{j1ort}\\
\langle P_{\bm{e}_j},P_{\bm{e}_k}\rangle&=0
\ \Longleftrightarrow \sum_{i=1}^n\eta_iu_{i\,j}u_{i\,k}=1,
\quad j\neq k,\quad j,k=1,\ldots,n.
\label{ijort}
\end{align}
These are obtained based on
\begin{equation}
P_{\bm{e}_j}(\bm{x})=1-\frac1N\sum_{i=1}^nx_iu_{i\,j},\quad \langle 1,x_i\rangle=\eta_iN,
\quad \langle x_j,x_k\rangle=\eta_j\eta_kN(N-1)+\eta_jN\delta_{j\,k}.
\label{deg1data}
\end{equation}
\end{rema}
%%%%%%%%%%%%%%%%%%%%%%%%%%%%%%%%%%%%%%%%%%%%%%%%%
The first $n$ conditions \eqref{j1ort} correspond to $\bigl(B^TD_1B\bigr)_{0\,j}=0$ as
\begin{align*}
0&=\bigl(B^TD_1B\bigr)_{0\,j}=\sum_{i=0}^n\bigl(B^T\bigr)_{0\,i}\bigl(D_1\bigr)_{i\,i}B_{i\,j}
=\bigl(B^T\bigr)_{0\,0}\bigl(D_1\bigr)_{0\,0}B_{0\,j}
+\sum_{i=1}^n\bigl(B^T\bigr)_{0\,i}\bigl(D_1\bigr)_{i\,i}B_{i\,j}\\
&=\eta_0+\sum_{i=1}^n\eta_ib_{i\,j}=\eta_0+\sum_{i=1}^n\eta_i(1-u_{i\,j})=\eta_0+\sum_{i=1}^n\eta_i
-\sum_{i=1}^n\eta_iu_{i\,j}=1-\sum_{i=1}^n\eta_iu_{i\,j}.
\end{align*}
Here $\eta_0+\sum_{i=1}^n\eta_i=1$ is used.
The second $n(n-1)/2$ conditions \eqref{ijort} correspond to $\bigl(B^TD_1B\bigr)_{j\,k}=0$, as
\begin{align*}
0&=\bigl(B^TD_1B\bigr)_{j\,k}=\sum_{i=0}^n\bigl(B^T\bigr)_{j\,i}\bigl(D_1\bigr)_{i\,i}B_{i\,k}
=\bigl(B^T\bigr)_{j\,0}\bigl(D_1\bigr)_{0\,0}B_{0\,k}
+\sum_{i=1}^n\bigl(B^T\bigr)_{j\,i}\bigl(D_1\bigr)_{i\,i}B_{i\,k}\\
& =\eta_0+\!\sum_{i=1}^nb_{i\,j}\eta_ib_{i\,k}=\eta_0+\sum_{i=1}^n\eta_i
-\sum_{i=1}^n\eta_iu_{i\,j}-\!\sum_{i=1}^n\eta_iu_{i\,k}+\sum_{i=1}^n\eta_iu_{i\,j}u_{i\,k}=\!-1\!+\sum_{i=1}^n\eta_iu_{i\,j}u_{i\,k},
\end{align*}
in which $\eta_0+\sum_{i=1}^n\eta_i=1$ and \eqref{j1ort} is used.
%%%%%%%%%%%%%%%%%%%%%%%%%%%%%%%%%%%%%%%%%%%%%%%
\begin{rema}
\label{GrRa}
In a slightly different context, in {\rm\cite{gr2}}, \eqref{j1ort} and its dual version are termed  the necessary condition
{\rm\cite{gr2}(3.3)}
of the orthogonality  of \eqref{Pm} and \eqref{ijort} is called the sufficient conditions {\rm\cite{gr2}(4.6)}
 of the orthogonality.
In {\rm\cite{Gri2}(25)} the system parameters $\{u_{i\,j}\}$ in the general hypergeometric function 
\eqref{Pm} are constrained by \eqref{j1ort} and \eqref{ijort}.
\end{rema}
%%%%%%%%%%%%%%%%%%%%%%%%%%%%%%%%%%%%%%%%%%%%%%%%
%%%%%%%%%%%%%%%%%%%%%%%%%%%%%%%%%%%%%%%%%%%%%%%
\begin{rema}
\label{challenge}
For the actual solutions of the difference equations, i.e., $\{\eta_j\}$ defined by \eqref{etadef} 
and $\{u_{i\,j}\}$  by \eqref{udef}, the conditions \eqref{j1ort} and \eqref{ijort} are actually satisfied.
It is rather easy to prove \eqref{j1ort} by the definitions \eqref{etadef}, \eqref{udef} and using \eqref{plambrel}.
It is an interesting challenge to prove \eqref{ijort} for general $n$ in a same way.
\end{rema}
%%%%%%%%%%%%%%%%%%%%%%%%%%%%%%%%%%%%%%%%%%%%%%%%

%%%%%%%%%%%%%%%%%%%%%%%%%%%%%%%%%%%%%
%
%         section3.6  Norm
%
%%%%%%%%%%%%%%%%%%%%%%%%%%%%%%%%%%%%%
\subsection{Norm}
\label{sec:norm}
According to {\bf Proposition \ref{mizuT}}, in order to find the norm of
general $\{P_{\bm{m}}(\bm{x})\}$,
one only has to evaluate the norm of the degree one
polynomials $\{P_{\bm{e}_j}(\bm{x})\}$, $j=1,\ldots,n$, 
\begin{equation}
\langle P_{\bm{e}_j},P_{\bm{e}_j}\rangle=\frac1N\Bigl(\sum_{i=1}^n\eta_iu_{i\,j}^2-1\Bigr)
\quad \Rightarrow\quad  \bar{p}_j=\Bigl(\sum_{i=1}^n\eta_iu_{i\,j}^2-1\Bigr)^{-1},\quad
\quad j=1,\ldots,n.
\tag{\ref{pddef}}
\end{equation}
Here, \eqref{deg1data} is used again. This leads to the norm formula \eqref{nKrorth}.
%%%%%%%%%%%%%%%%%%%%%%%%%%%%%%%%%%%%%
%
%         section3.7  Dual Polynomials
%
%%%%%%%%%%%%%%%%%%%%%%%%%%%%%%%%%%%%%
\subsection{Dual Polynomials}
\label{sec:dual}
In terms of the norm formula \eqref{nKrorth}, the following  set of orthonormal vectors on $\mathcal{X}$ are
defined
\begin{align}
&\quad \sum_{\bm{x}\in\mathcal{X}}\hat{\phi}_{\bm{m}}(\bm{x})\hat{\phi}_{\bm{m}'}(\bm{x})
=\delta_{\bm{m}\,\bm{m}'}, \qquad \quad \bm{m},\bm{m}'\in\mathcal{X},
\label{ortrel}\\
\hat{\phi}_{\bm{m}}(\bm{x})&\eqdef \sqrt{W(\eta;\bm{x})}P_{\bm{m}}(\bm{x})\sqrt{\bar{W}(\bar{p};\bm{m})},
\qquad \bm{x},\bm{m}\in\mathcal{X},
\label{ortphiderf}\\
&\bar{W}(\bar{p};\bm{m})\eqdef \binom{N}{\bm{m}}(\bar{\bm{p}})^{\bm{m}},\quad
\sum_{\bm{m}\in\mathcal{X}}\bar{W}(\bar{p};\bm{m})=\Bigl(1+\sum_{j=1}^n\bar{p}_j\Bigr)^N.
\label{baretaW}
\end{align}
They define an orthogonal matrix $\mathcal{T}$ on $\mathcal{X}$,
\begin{align*}
&\hspace{4cm} \mathcal{T}_{\bm{x}\,\bm{m}}\eqdef \hat{\phi}_{\bm{m}}(\bm{x}),\\
&\bigl(\mathcal{T}^T\mathcal{T}\bigr)_{\bm{m}\,\bm{m}'}=
\sum_{\bm{x}\in\mathcal{X}}\bigl(\mathcal{T}^T\bigr)_{\bm{m}\,\bm{x}}\mathcal{T}_{\bm{x}\,\bm{m}'}
=\sum_{\bm{x}\in\mathcal{X}}\hat{\phi}_{\bm{m}}(\bm{x})\hat{\phi}_{\bm{m}'}(\bm{x})
=\delta_{\bm{m}\,\bm{m}'},\n
&\Rightarrow \delta_{\bm{x}\,\bm{y}}=\bigl(\mathcal{T}\mathcal{T}^T)_{\bm{x}\,\bm{y}}
=\sum_{\bm{m}\in\mathcal{X}}\mathcal{T}_{\bm{x}\,\bm{m}}\bigl(\mathcal{T}^T\bigr)_{\bm{m}\,\bm{y}}
=\sum_{\bm{m}\in\mathcal{X}}\hat{\phi}_{\bm{m}}(\bm{x})\hat{\phi}_{\bm{m}}(\bm{y}).
\end{align*}
This means that $\hat{\phi}_{\bm{m}}(\bm{x})$ defines dual polynomials in $\bm{m}$ indexed by $\bm{x}$,
%%%%%%%%%%%%%%%%%%%%%%%%%%%%%%%%%%%%%%%%%%%%%%%%%%%%
\begin{prop}
\label{dualQ}
The dual polynomials of $\{P_{\bm{m}}(\bm{x})\}$ to be denoted by $\{Q_{\bm{x}}(\bm{m})\}$ are defined by the same
formula as $P_{\bm{m}}(\bm{x})$ \eqref{Pm},
\begin{gather}
\label{Qm}
Q_{\bm{x}}(\bm{m})
\eqdef \sum_{\substack{\sum_{i,j}c_{ij}\leq N\\
(c_{ij})\in M_{n}({\mathbb N_{0}})}}
\frac{\prod\limits_{i=1}^{n}(-x_{i})_{\sum\limits_{j=1}^{n}c_{ij}}
\prod\limits_{j=1}^{n}(-m_{j})_{\sum\limits_{i=1}^{n}c_{ij}}}
{(-N)_{\sum_{i,j}c_{ij}}} \; \frac{\prod(u_{ij})^{c_{ij}}}{\prod c_{ij}!}.
\end{gather}
They are orthogonal with respect to the dual multinomial distribution
\begin{align}
&W(\bar{\eta};\bm{m})\eqdef \binom{N}{\bm{m}}\bigl(\bar{\eta}_0\bigr)^{m_0}
\bigl(\bar{\bm{\eta}}\bigr)^{\bm{m}},
\quad \sum_{\bm{m}\in\mathcal{X}}W(\bar{\eta};\bm{m})=1,
\label{Wddef}\\
&\sum_{\bm{m}\in\mathcal{X}}W(\bar{\eta};\bm{m})Q_{\bm{x}}(\bm{m})Q_{\bm{y}}(\bm{m})
=\frac{\delta_{\bm{x}\,\bm{y}}}{W(\eta;\bm{x})\bigl(\bar{\eta}_0\bigr)^{-N}},\qquad \bm{x},\bm{y}\in\mathcal{X},
\label{Qmorth}
\end{align}
in which
%\begin{equation}
%\eta_0^d\eqdef\bigl(1+\sum_{j=1}^n\bar{\eta}_j\bigr)^{-1},\quad
%\eta_j^d\eqdef\eta_0^d\bar{\eta}_j,\quad j=1,\ldots,n.
%\end{equation}
\begin{equation}
\bar{\eta_i}\eqdef\frac{\bar{p}_i}{1+\sum_{j=1}^n\bar{p}_j},\quad 
\bar{\eta}_0\eqdef\frac1{1+\sum_{i=1}^n\bar{p}_i},
\quad \sum_{i=0}^n\bar{\eta}_i=1, \quad (\bar{\bm{\eta}})^{\bm{m}}\eqdef\prod_{i=1}^n(\bar{\eta}_i)^{m_i}.
\end{equation}

\end{prop}
%%%%%%%%%%%%%%%%%%%%%%%%%%%%%%%%%%%%%%%%%%%%%%%%%%%%%%%%%%
The orthogonalities of $\{Q_{\bm{x}}(\bm{m})\}$ are expressed, corresponding to \eqref{j1ort}, \eqref{ijort}, as
\begin{align}
 \sum_{i=1}^n\bar{\eta}_iu_{j\,i}=1,
\qquad j=1,\ldots,n,
\label{dj1ort}\\
 \sum_{i=1}^n\bar{\eta}_iu_{j\,i}u_{k\,i}=1,
\quad j\neq k,\quad j,k=1,\ldots,n.
\label{dijort}
\end{align}
The above dual orthogonalities are the consequences of the original orthogonalities \eqref{j1ort}, \eqref{ijort},
not the conditions.
%%%%%%%%%%%%%%%%%%%%%%%%%%%%%%%%%%%%%%%%%%%%%%%%%%%%%%
\begin{rema}
\label{dualBD}
It is obvious that the dual polynomials $\{Q_{\bm{x}}(\bm{m})\}$ are also governed by
dual birth and death process with the BD rates
\begin{equation}
B_j^d(\bm{m})=\bigl(N-|m|),\quad D_j^d(\bm{m})=(p_j^d)^{-1}m_j,\quad j=1,\ldots,n.
\label{dBD}
\end{equation}
with certain dual probability parameters $\{p_j^d\}$, $j=1,\ldots,n$.
In sharp contrast to the single variable BD cases, the relationship connecting  the dual parameters $\{p_j^d\}$ 
with  the original BD parameters $\{p_j\}$ seems highly nontrivial.
\end{rema}
%%%%%%%%%%%%%%%%%%%%%%%%%%
%%%%%%%%%%%%%%%%%%%%%%%%%%%%%%%%%%%%%
%
%         section3.8  Exceptional Cases
%
%%%%%%%%%%%%%%%%%%%%%%%%%%%%%%%%%%%%%
\subsection{Exceptional Cases}
\label{sec:excep}
So far the parameter values are assumed to be generic. 
But obviously at certain parameter settings, 
the above hypergeometric formula \eqref{Pm} for the polynomials $\{P_{\bm{m}}(\bm{x})\}$ could go wrong. 
By construction, the eigenvalues $\{\lambda_i\}$ are positive 
and $\{p_i\}$ are also positive. 
Therefore, if the situation $\lambda_j=p_i^{-1}$  happens at some parameter setting, it leads to the
breakdown of the generic theory  as 
$u_{i\,j}=\frac{\lambda_j}{\lambda_j-p_i^{-1}}$  \eqref{udef} is ill-defined.

%%%%%%%%%%%%%%%%%%%%%%%%%%%%%%%%%%%%%%%%%%%%%
\subsubsection{$n=2$ Case}
\label{sec:n2exc}
The situation is most clearly seen when $n=2$. In this case the two eigenvalues are
the roots of the characteristic equation
\begin{align*}
&\lambda^2-(2+p_1^{-1}+p_2^{-1})\lambda+(1+p_1^{-1})(1+p_2^{-1})-1=0,\n
\lambda_1&=\frac12(2+p_1^{-1}+p_2^{-1}-\Delta),\quad \lambda_2=\frac12(2+p_1^{-1}+p_2^{-1}+\Delta),\\
& \Delta^2=4+(1/p_1-1/p_2)^2.
\end{align*}
When $p_1=p_2=p$, the eigenvalues are rational,
\begin{equation*}
p_1=p_2=p\ \Longrightarrow \lambda_1=1/p,\quad \lambda_2=2+1/p,
\end{equation*}
and the singular situation occurs, $\lambda_1-1/p=0$. That is, the general formula \eqref{Pm}  fails.
In this case, the degree 1 solution of $\widetilde{\mathcal H}$ \eqref{HthnKr}
corresponding to the eigenvalue $\lambda_1$ is
\begin{equation}
P_{\bm{e}_1}(\bm{x})=const\times\bigl(x_1-x_2).
\label{non1sol}
\end{equation}
That is the constant part is vanishing and the assumption that degree one solutions have unit 
constant part \eqref{deg1form} simply fails. The situation is similar for general $n$ as 
stated by the following 
%%%%%%%%%%%%%%%%%%%%%%%%%%%%%%%%%%%%%%%%%%%%%
\begin{theo}
\label{theo:sing}
When some of the parameters $\{p_i\}$ coincide, the hypergeometric formula \eqref{Pm}
for the $n$-variate Krawtchouk polynomials does not apply.
But the solutions of the difference equations $\widetilde{\mathcal H}$ \eqref{HthnKr}
still constitute the $n$-variate orthogonal polynomials.
In other words, the hypergeometric formula \eqref{Pm} requires all distinct probability parameters
$\{p_i\}$.
\end{theo}
This is rather easy to see. If $p_j=p_k=p$, the matrix $p^{-1}I_n-F(p)$ \eqref{chareq} 
has the $j$-th and $k$-th column $-(1,1,\ldots,1)^T$,
thus the characteristic polynomial $\mathcal{F}(\lambda)$ vanishes at $\lambda=p^{-1}$.
When  there exist $k$ identical $p_i$'s, $\mathcal{F}(\lambda)$ has a factor $(\lambda-p_i^{-1})^{k-1}$.
%%%%%%%%%%%%%%%%%%%%%%%%%%%%%%%%%%%%%%%%%%%%%%%
%%%%%%%%%%%%%%%%%%%%%%%%%%%%%%%%%%%%%%%%%%%%%%%
\begin{theo}
\label{theo:distinct}{\bf Distinct parameters $\{p_j\}$ are necessary}\\
All the parameters $\{p_j\}$ must be distinct for the hypergeometric formula \eqref{Pm}
for the $n$-variate Krawtchouk polynomials to hold.
\end{theo}
%%%%%%%%%%%%%%%%%%%%%%%%%%%%%%%%%%%%%%%%%%%%%%%
\begin{rema}
\label{fullform}
It is a big challenge to derive a general formula of $n$-variate Krawtchouk polynomials 
including all these exceptional cases.
\end{rema}
%%%%%%%%%%%%%%%%%%%%%%%%%%%%%%%%%%%%%%%%%%%%%
%%%%%%%%%%%%%%%%%%%%%%%%%%%%%%%%%%%%%%%%%%%%%
%\begin{prop}{\bf $n=2$  Rational Case}
%\label{prop:n2rat}
%As is clear from \eqref{2eig2}, for $q_1=q$, $q_2=q+2(p_1-p_2)$, the eigenvalues are rational
%$\lambda_1=q+p_1-p_2$, $\lambda_2=q+2p_1$ and the 2-variate Krawtchouk polynomials have 
%rational coefficients only
%\begin{equation}
%u_{1\,1}=\frac{q+p_1-p_2}{p_1-p_2},\quad u_{1\,2}=\frac{q+2p_1}{2p_1},
%\quad u_{2\,1}=\frac{q+p_1-p_2}{p_2-p_1},\quad u_{2\,2}=\frac{q+2p_1}{2p_2}.
%\end{equation}
%\end{prop}

%%%%%%%%%%%%%%%%%%%%%%%%%%%%%%%%%%%%%
%
%         section4  Rahman Polynomials
%
%%%%%%%%%%%%%%%%%%%%%%%%%%%%%%%%%%%%%
\section{Rahman Polynomials}
\label{sec:rahman}

It is obvious that the Rahman polynomials \cite{HR, gr1, gr2, gr3, IT, I} and 
the multivariate Krawtchouk polynomials presented in this paper share the bulk of the basic structure.
It is rather intricate to demonstrate the actual relationship for the general $n$-variate cases,
especially due to the very special choice of parameters in the multivariate Rahman polynomials.
However, detailed comparison is possible for the bivariate Rahman polynomials due to  following 
%%%%%%%%%%%%%%%%%%%%%%%%%%%%%%%%%%%%%%
\begin{theo}
\label{theo:grurahman}
{\bf Gr\"unbaum and Rahman \cite{gr2}} showed that the bivariate Rahman polynomials 
\begin{equation*}
P_{m,n}(x,y)=
{\underset{0 \le i+j+k+l \le N}{\sum_i\sum_j\sum_k\sum_l}} \frac {(-m)_{i+j}(-n)_{k+l}(-x)_{i+k}(-y)_{j+l}}{i!j!k!l!(-N)_{i+j+k+l}} t^iu^jv^kw^l, \tag{\cite{gr2}.1.2}
\end{equation*}
satisfy the orthogonality {\rm\cite{gr2}(1.8)} and the dual orthogonality {\rm\cite{gr2}(2.4)} relations and a 
5-term recurrence relation,
\begin{gather*}
(N-m_1-m_2) \left\{ \frac {p_1p_3(p_2+p_4)(p_1+p_2+p_3+p_4)}{(p_1+p_3)(p_1p_4-p_2p_3)} (P_{m_1+1,m_2}(x_1,x_2) - P_{m_1,m_2}(x_1,x_2))\right. \\
\qquad \left. {} - \frac {p_2p_4(p_1+p_3)(p_1+p_2+p_3+p_4)}{(p_2+p_4)(p_1p_4 - p_2p_3)} (P_{m_1,m_2+1}(x_1,x_2) - P_{m_1,m_2}(x_1,x_2))\right\} \\
 \qquad {} + m_1 \frac {p_1p_4-p_2p_3}{p_1+p_3} (P_{m_1-1,m_2}(x_1,x_2) - P_{m_1,m_2}(x_1,x_2)) \\
 \qquad {} - m_2 \frac {p_1p_4-p_2p_3}{p_2+p_4} (P_{m_1,m_2-1}(x_1,x_2) - P_{m_1,m_2}(x_1,x_2)) \\
 \qquad {} = ((p_1+p_2)x_1 - (p_3+p_4)x_2) P_{m_1,m_2}(x_1,x_2).
 \tag{\cite{gr2}.1.9}
\end{gather*}
which is, in fact, the difference equation governing the dual birth and death process.
\end{theo}
By reversing the logic, I derive these results starting from the dual BD difference equation.
%%%%%%%%%%%%%%%%%%%%%%%%%%%%%%%%%%%%%%
\begin{theo}
\label{theo:2rahman}
The bivariate Rahman polynomials  in {\rm\cite{gr1,gr2}} are dual birth and death 
polynomials for a very special choice of the birth and death parameters.
They are denoted tentatively by $\{Q_{\bm{x}}(\bm{m})\}$, as they are dual polynomials,
\begin{align}
%&P_{\bm{m}}(\bm{x})=\sum_{i}\sum_{j}\sum_{k}\sum_{l}
%\frac{(-m_1)_{i+j}(-m_2)_{k+l}(-x_1)_{i+k}(-x_2)_{j+l}}
%{i!j!k!l!(-N)_{i+j+k+l}}\,t^iu^jv^kw^l,\n
%
&\widetilde{\mathcal  H}^dQ_{\bm{x}}(\bm{m})
=\Bigl(-(p_1+p_2)x_1+(p_3+p_4)x_2\Bigr)Q_{\bm{x}}(\bm{m}),
\end{align}
in which the dual difference operator is
\begin{align}
&\widetilde{\mathcal  H}^d\eqdef(N-m_1-m_2)
\bigl(p_1^d(1-e^{\partial_1^d})+p_2^d(1-e^{\partial_2^d})\bigr)
\!+q_1^dm_1(1-e^{-\partial_1^d})+q_2^dm_2(1-e^{-\partial_2^d}),
\label{Hthddef}\\
& \hspace{4cm} \partial_1^d\eqdef\frac{\partial}{\partial m_1},\quad 
 \partial_2^d\eqdef\frac{\partial}{\partial m_2},\n
&\qquad p_1^d\eqdef \frac{p_1p_3(p_2+p_4)(p_1+p_2+p_3+p_4)}{(p_1+p_3)(p_1p_4-p_2p_3)},
\qquad \ \ q_1^d\eqdef \frac{p_1p_4-p_2p_3}{p_1+p_3},
\label{dpardef1}\\
&\qquad p_2^d\eqdef-\frac{p_2p_4(p_1+p_3)(p_1+p_2+p_3+p_4)}{(p_2+p_4)(p_1p_4-p_2p_3)},
\qquad q_2^d\eqdef -\frac{p_1p_4-p_2p_3}{p_2+p_4}.
\label{dpardef2}
\end{align}
The dual parameters $p_1^d,q_1^d,p_2^d,q_2^d$ are taken from  above {\rm\cite{gr2}(1.9)}.
The operator $-\widetilde{\mathcal  H}^d$ corresponds to l.h.s. of {\rm\cite{gr2}(1.9)} 
and to the matrix $\mathcal{B}$ in \S6 of {\rm\cite{gr1}}.
\end{theo}
%%%%%%%%%%%%%%%%%%%%%%%%%%%%%%%%%%%%%%%%%%%%%%
This produces all the other system parameters as shown below. 
First, the dual probabilities which are listed as $\bar{\eta}_1$ and $\bar{\eta}_2$ in \cite{HR} are
\begin{align*}
\eta_1^d&=\frac{p_1^d/q_1^d}{1+\sum_{i=1}^2p_i^d/q_i^d}
=\frac{p_1p_3(p_1+p_2+p_3+p_4)}{(p_1+p_2)(p_1+p_3)(p_3+p_4)}=\bar{\eta}_1,
\tag{\cite{HR}.2.26}\\
\eta_2^d&=\frac{p_2^d/q_2^d}{1+\sum_{i=1}^2p_i^d/q_i^d}
=\frac{p_2p_4(p_1+p_2+p_3+p_4)}{(p_1+p_2)(p_2+p_4)(p_3+p_4)}=\bar{\eta}_2.
\tag{\cite{HR}.2.27}
\end{align*}
The degree one eigenpolynomials of $\widetilde{\mathcal  H}^d$ \eqref{Hthddef} 
are determined in a similar way to
\eqref{deg1form},
\begin{align*}
&Q_{|x|=1}(\bm{m})=1+\sum_{i=1}^2b_im_i,
\quad
\widetilde{\mathcal H}^dQ_{|x|=1}(\bm{m})=\lambda^d Q_{|x|=1}(\bm{m}),\\
%\label{deg1form}\\
%
&\Rightarrow -(N-m_1-m_2)\sum_{i=1}^2p_i^db_i+\sum_{i=1}^2q_i^db_im_i
=\lambda^d(1+\sum_{i=1}^2b_im_i).
\end{align*}
The  characteristic equation corresponding to \eqref{chareq} reads
\begin{align}
&(\lambda^d)^2-(p_1^d+q_1^d+p_2^d+q_2^d)\lambda^d+(p_1^d+q_1^d)(p_2^d+q_2^d)-p_1^dp_2^d
=(\lambda^d+p_1+p_2)(\lambda^d-p_3-p_4),\n
&\quad \Longrightarrow \lambda_1^d=-(p_1+p_2),\quad \lambda_2^d=p_3+p_4.
\label{dlamb}
\end{align}
Corresponding to \eqref{e_jsol}, the eigenpolynomial related to
$\bm{x}=\bm{e}_j$ are
\begin{align*}
&Q_{\bm{e}_j}(\bm{m})=1-\frac1N\sum_{i=1}^2m_iu_{j\,i}^d,\quad 
u_{j\,i}^d\eqdef \frac{\lambda_j^d}{\lambda_j^d-q_i^d},\quad
i,j=1,2,\\
&u_{1\,1}^d=\frac{(p_1+p_2)(p_1+p_3)}{p_1(p_1+p_2+p_3+p_4)}\equiv t,\quad
u_{1\,2}^d=\frac{(p_1+p_2)(p_2+p_4)}{p_1(p_1+p_2+p_3+p_4)}\equiv v,
\tag{\cite{gr2}.1.3}\\
&u_{2\,1}^d=\frac{(p_1+p_3)(p_3+p_4)}{p_1(p_1+p_2+p_3+p_4)}\equiv u,\quad
u_{2\,2}^d=\frac{(p_4+p_2)(p_4+p_3)}{p_1(p_1+p_2+p_3+p_4)}\equiv w.
\tag{\cite{gr2}.1.3}
\end{align*}
The parameters $\{u_{i\,j}^d\}$ determine the entire polynomials as \eqref{Qm}
and the parameters $t, u, v, w$ determine the above Rahman polynomials  \cite{gr2}(1.2).
%\begin{equation*}
%P_{m,n}(x,y)=\sum_{i}\sum_{j}\sum_{k}\sum_{l}
%\frac{(-m)_{i+j}(-n)_{k+l}(-x)_{i+k}(-y)_{j+l}}
%{i!j!k!l!(-N)_{i+j+k+l}}\,t^iu^jv^kw^l.
%\tag{\cite{HR}.(2.28)}
%\end{equation*}
Of course, the dual orthogonality conditions corresponding to \eqref{j1ort} and \eqref{ijort}
are satisfied,
\begin{equation*}
\sum_{i=1}^2\eta_i^du_{j\,i}^d=1, \quad j=1,2,\qquad \sum_{i=1}^2\eta_i^du_{1\,i}^du_{2\,i}^d=1.
\tag{\cite{gr2}.2.4}
\end{equation*}
By evaluating the norm of $Q_{\bm{e}_j}(\bm{m})$ in a similar way to \eqref{pddef}, one obtains
the probabilities of the original $\bm{x}$ system, $\eta_0$, $\eta_1$ and $\eta_2$,
\begin{align*}
\bar{\eta}_j^d&\eqdef\Bigl(\sum_{i=1}^2\eta_iu_{j\,i}^2-1\Bigr)^{-1},\quad
\bar{\eta}_1^d=\frac{p_1p_2(p_3+p_4)(p_1+p_2+p_3+p_4)}{(p_1p_4-p_2p_3)^2},\\
&\hspace{4.4cm}\bar{\eta}_2^d=\frac{p_3p_4(p_1+p_2)(p_1+p_2+p_3+p_4)}{(p_1p_4-p_2p_3)^2},\\
&\eta_0\eqdef\bigl(1+\bar{\eta}_1^d+\bar{\eta}_2^d)^{-1}=\frac{(p_1p_4-p_2p_3)^2}{(p_1+p_2)(p_1+p_3)(p_2+p_4)(p_3+p_4)}, \tag{\cite{HR}.2.22)}\\
&\eta_1\eqdef\bar{\eta}_1^d\eta_0=\frac{p_1p_2(p_1+p_2+p_3+p_4)}{(p_1+p_2)(p_1+p_3)(p_2+p_4)},
\tag{\cite{gr2}.1.4}\\
&\eta_2\eqdef\bar{\eta}_1^d\eta_0=\frac{p_3p_4(p_1+p_2+p_3+p_4)}{(p_1+p_3)(p_2+p_4)(p_3+p_4)}.
\tag{\cite{gr2}.1.4}
\end{align*}
The original $\bm{x}$ orthogonality conditions are satisfied as \eqref{j1ort} and \eqref{ijort}
\begin{equation*}
\sum_{i=1}^2\eta_iu_{i\,j}^d=1, \quad j=1,2,\qquad \sum_{i=1}^2\eta_iu_{i\,1}^du_{i\,2}^d=1.
\tag{\cite{gr2}.1.8}
\end{equation*}
In \cite{gr2} the 5-term recurrence relation \cite{gr2}(1.9) was proved by direct calculation, instead of using the
generating function as in {\bf Proposition \ref{lineareig}}.
%%%%%%%%%%%%%%%%%%%%%%%%%%%%%%%%%%%%%%%%%%%%%%%%%%%
\begin{rema}
\label{rpara}
The  original parameters $p_1$, $p_2$, $p_3$ and $p_4$ of the bivariate Rahman polynomials
are very special in the sense that the eigenvalues and the parameters describing the polynomials,
$u_{i\,j}^d$ or $t$, $u$, $v$ and $w$ are all rational functions of them. 
As shown in \S\ref{sec:nKr}, the default parametrisation of the multivariate BD process inevitably 
involves the irrational system parameters $u_{i\,j}$ \eqref{udef}. 
%For the bivariate case, a simple rational parametrisation contains only 3 degrees of freedom as shown in 
%{\bf Proposition \ref{prop:n2rat}}. 
It seems rather difficult to find the parametrisation 
of multivariate Rahman polynomials in which the system parameters $u_{i\,j}^d$, $i,j=1,\ldots,n$
are all rational, as $u_{i\,j}^d$ involves the eigenvalues of an $n\times n$ matrix eigenvalue problem.
The rational parametrisation of the bivariate Rahman polynomials seems to be the fortuitous outcome of 
of the 9-$j$ symbol's origin {\rm\cite{HR}}. 
\end{rema}
%%%%%%%%%%%%%%%%%%%%%%%%%%%%%%%%
%%%%%%%%%%%%%%%%%%%%%%%%%%%%%%%%%%%%%%%%%%%%%%%%%%%
\begin{rema}
\label{rpara2}
The ranges of the original parameters $p_1$, $p_2$, $p_3$ and $p_4$ of the bivariate Rahman polynomials
are unclear due to their limiting process origin.
For all positive range $p_i>0$, the corresponding trinomial probability $\eta_i$ {\rm\cite{gr2}(1.4)}
and the dual probability $\bar{\eta}_i$ {\rm\cite{HR}(2.26),(2.27)} are all positive.
However, some of the corresponding birth and death parameters  \eqref{dpardef1}, \eqref{dpardef2}  
and one of the eigenvalues \eqref{dlamb} are negative,
\begin{equation*}
p_1^dp_2^d<0,\quad q_1^dq_2^d<0,\quad -(p_1+p_2)(p_3+p_4)<0.
\end{equation*}
This means that the corresponding BD process is explosive, that is it does not tend to the
stationary probability distribution.
It is interesting to see if this causes any trouble or not in the actual applications of the bivariate
Rahman polynomials.
\end{rema}

%%%%%%%%%%%%%%%%%%%%%%%%%%%%%%%%%%%%%
%
%         section5 Comments
%
%%%%%%%%%%%%%%%%%%%%%%%%%%%%%%%%%%%%
\section{Comments}
\label{sec:comm}
Physicists knew that the Schr\"odinger equations and the Fokker-Planck equations 
are related by similarity transformations \cite{risken}.
That is, the eigenvalue problem of a self-adjoint operator is related to an equation governing the diffusion processes.
The same mathematical structure, a discretised version, exists between the BD processes and the difference
equations governing the orthogonal polynomials.

Once the explicit forms of the multivariate orthogonal polynomials with many free parameters are available, 
a lot of interesting questions pop up and demand answers. 
The situation looks much more complex and interesting than 
the Calogero-Moser-Sutherland systems which have only a few free parameters.
Here I name only a few.
\begin{itemize}
\item How the oscillation theorem, if any, can be formulated? 
How are the numbers of the positive, negative regions or the boundaries related with the degree $\bm{m}$, $|m|$ and/or
the eigenvalues $\lambda_i$, $\mathcal{E}(\bm{m})$?
In this connection, is it meaningful and useful to order the eigenvalues $\lambda_1<\lambda_2<\cdots<\lambda_n$?
\item Level crossing? What happens, by parameter changes, when the eigenvalues cross, $\lambda_i<\lambda_j\to \lambda_j<\lambda_i$?
\item What happens when two death rates cross $q_i<q_j\to q_j<q_i$?
\end{itemize}
Of course it is a great challenge to enlarge the list of multivariate discrete orthogonal polynomials which are
the generalisation of the Askey scheme polynomials.

The fact that the divariate Rahman polynomials have four independent parameters on top of $N$
led me to suppose that the general multivariate Krawtchouk polynomials have $2n+1$ parameters.

After the first version of this work was published, Plamen Iliev informed me that he had arrived at the same polynomials
with $2n+1$ parameters \cite{I23} several weeks earlier.

%%%%%%%%%%%%%%%%%%%%%%%%%%%%%%%%%%%%%%%%%%%%%%%%%%%%%%%%%%%%%%%
%                                                             %
%  Acknowledgments                                            %
%                                                             %
%%%%%%%%%%%%%%%%%%%%%%%%%%%%%%%%%%%%%%%%%%%%%%%%%%%%%%%%%%%%%%%
\section*{Acknowledgements}
RS thanks Mourad Ismail for inducing him to explore the mysterious maze of birth and death processes.
He also thanks Plamen Iliev for the information of his recent work.
%%%%%%%%%%%%%%%%%%%%%%%%%%%%%%%%%%%%%%%%%%%%%%%%%%%%%
%                                                             %
%  References                                                 %
%                                                             %
%%%%%%%%%%%%%%%%%%%%%%%%%%%%%%%%%%%%%%%%%%%%%%%%%%%%

%\goodbreak

\begin{thebibliography}{99}

\bibitem{AK}
Aomoto K., Kita M.,
{\it Hypergeometric funtions}, Springer, Berlin, 1994 (in Japanese).


\bibitem{askey}
Andrews G.E.,\, Askey R.\, and Roy R.,
{\it Special Functions},
Encyclopedia of mathematics and its applications,
Cambridge Univ. Press, Cambridge, (1999).


\bibitem{coo-hoa-rah77}
Cooper R.\,D., Hoare M.\,R. and Rahman M.,
``Stochastic Processes and Special Functions:
On the Probabilistic Origin of Some Positive Kernels Associated
with Classical Orthogonal Polynomials,"
J. Math. Anal. Appl. {\bf 61} (1977) 262-291.


\bibitem{diaconis13}
Diaconis P.\, and Grif\/f\/iths R.C.,
``An introduction to multivariate Krawtchouk polynomials and their applications,"
J. Stat. Planning and Inference {\bf 154} (2014) 39-53,
{\tt arXiv:\hspace{0pt}1309.0112[math.PR]}.

\bibitem{feller}
Feller W., {\it An Introduction to Probability Theory and its Applications,
I}, (2nd ed.), Wiley, New York, (1966).

\bibitem{gasper}
Gasper G.\, and Rahman M.,
{\it Basic hypergeometric series\/}, 2nd ed.
Encyclopedia of mathematics and its applications, Cambridge, (2004).


 \bibitem{gelfand}
 Gelfand I.\,M., ``General theory of hypergeometric functions,"  Sov. Math. Dokl. {\bf 33} (1986) 573-577.

\bibitem{genest}
Genest V.\,X., \ Miki H.,\, Vinet L. and Zhedanov A.,
``Spin lattices, state transfer and bivariate Krawtchouk polynomials,"
J. Phys. A: Math. Theor. {\bf 46} (2013) 505203,
{\tt arXiv:\hspace{0pt}1410.4703[math-ph]}.

\bibitem{Gri1}
Grif\/f\/iths R.\,C.,
``Orthogonal polynomials on the multinomial distribution,"
%\href{http://dx.doi.org/10.1111/j.1467-842X.1971.tb01239.x}
{{\it Austral.~J. Statist.}} {\bf 13} (1971)   27--35.

\bibitem{Gri2}
Grif\/f\/iths R.\,C.,
``Multivariate Krawtchouk Polynomials and Composition Birth and Death Processes,"
Symmetry  {\bf 8} (2016) 33 19pp. {\tt arXiv:1603.00196[math.PR]}.


\bibitem{gr1}
 Gr\"{u}nbaum F.\,A.,
``The Rahman polynomials are bispectral,"
 SIGMA {\bf 3} (2007) 065, 11pp,
 {\tt arXiv:0705.0468[math.CA]}.
 
 
\bibitem{gr2}
 Gr\"{u}nbaum F.\,A. and Rahman M.,
``On a family of 2-variable orthogonal Krawtchouk polynomials,"
%\href{http://dx.doi.org/10.3842/SIGMA.2010.090}
SIGMA {\bf 6} (2010) 090, 12~pages,
%\href{http://arxiv.org/abs/1007.4327}
{\tt arXiv:1007.4327[math.CA]}.

\bibitem{gr3}
 Gr\"{u}nbaum F.\,A. and  Rahman M.,
 ``A System of multivariable Krawtchouk polynomials
and a probabilistic application,"
SIGMA {\bf 7} (2011) 119, 17pp,
{\tt arXiv:1106.1835\hspace{0pt}[math.PR]}.

\bibitem{HR}
Hoare M.\,R. and  Rahman M.,
``A probabilistic origin for a new class of bivariate polynomials,"
%\href{http://dx.doi.org/10.3842/SIGMA.2008.089}
SIGMA {\bf 4} (2008)  089, 18~pages,
%\href{http://arxiv.org/abs/0812.3879}
{arXiv:0812.3879[math.CA]}.

\bibitem{illxu}
Iliev P. and Xu Y.,
``Discrete orthogonal polynomials and difference equations of several variables,"
 Adv. Math. { \bf 212} (2007) 1-36, {\tt arXiv:math.CA/0508039}.

\bibitem{IT}
Iliev P. and  Terwilliger P.,
``The Rahman polynomials and the Lie algebra $sl_3(C)$,"
Trans. Amer. Math. Soc. {\bf 364} (2012) 4225--4238,
%\href{http://arxiv.org/abs/1006.5062}
{\tt arXiv:1006.5062[math.RT]}.

\bibitem{I}
Iliev P.,
``A Lie theoretic interpretation of multivariate hypergeometric polynomials,"
Compositio Math. {\bf 148} (2012) 991-1002,
%\href{http://arxiv.org/abs/1101.1683}
{\tt arXiv:1101.1683[math.RT]}.


\bibitem{I23}
Iliev P.,
``Gaudin model for the multinomial distribution,"
%\href{http://arxiv.org/abs/2303.08206}
{\tt arXiv:2303.08206[math-\hspace{0pt}ph]}.



\bibitem{ismail}
Ismail M.\,E.\,H.,
{\it Classical and quantum orthogonal polynomials in one variable\/},
Encyclopedia of mathematics and its applications, Cambridge, (2005).



\bibitem{KarMcG}
Karlin S.  and \,McGregor J.\,L., ``The differential equations
of birth-and-death processes," Trans. Amer. Math. Soc. {\bf 85} (1957)
489-546;
``Linear growth, birth-and-death processes," 
J.  Math. Mech. {\bf 7} (1958) 643--662;
``Ehrenfest urn models," 
J.  Appl. Prob. {\bf 19} (1965) 477--487.


\bibitem{khare}
Khare K. and Zhou  H., ``Rates of convergence of some multivariate Markov chains 
with polynomial eigenfucntions," Ann. Appl. Probab. {\bf 19}  (2009) 737-777,
{\tt arXiv:\hspace{0pt}0906.4242[math.PR]}.

\bibitem{koeswart} 
Koekoek R.,  Lesky P.\,A.  and Swarttouw R.\,F.,
{\it Hypergeometric orthogonal polynomials and their $q$-analogues,\/}
Springer Monographs in Mathematics, 
Springer-Verlag, Berlin, (2010).
%

\bibitem{mizu1}
 Mizukawa H.,
 ``Zonal spherical functions on the  complex reflection groups and $(m+1,n+1)$-hypergeometric functions,"
%\href{http://dx.doi.org/10.1016/S0001-8708(03)00092-6}
{\it Adv. Math.} {\bf 184} (2004) 1--17.

\bibitem{mizu}
 Mizukawa H.,
``Orthogonal relations for multivariate Krawtchouk polynomials,"
SIGMA {\bf 7} (2011) 017 5pp,
{\tt arXiv:1009.1203[math.CO]}.



 \bibitem{mt}
Mizukawa H. and  Tanaka H.,
 $(n+1,m+1)$-hypergeometric functions associated to character algebras,
%\href{http://dx.doi.org/10.1090/S0002-9939-04-07399-X}
{{\it Proc. Amer. Math. Soc.}} {\bf 132} (2004) 2613--2618.

\bibitem{risken}
Risken H., {\it The Fokker-Planck Equation}, second ed., Springer-Verlag, Berlin, (1996).

\bibitem{os12}
Odake S. and Sasaki R.,
``Orthogonal Polynomials from Hermitian Matrices,"
J. Math. Phys. {\bf 49} (2008) 053503 (43 pp),
{\tt arXiv:0712.4106[math.CA]}.


\bibitem{os34}
Odake S. and Sasaki R.,
``Orthogonal Polynomials from Hermitian Matrices II,"
J. Math. Phys. {\bf 59} (2018) 013504 (42pp)
{\tt arXiv:1604.00714[math.CA]}. 

\bibitem{bdsol}
Sasaki R.,
``Exactly Solvable Birth and Death Processes,"
J. Math. Phys. {\bf 50} (2009) 103509 (18 pp),
 {\tt arXiv:0903.3097[math-ph]}.

\bibitem{tra}
Tratnik M.V., ``Some multivariable orthogonal polynomials of the Askey tableau-discrete families,"
 J. Math. Phys. {\bf 32} (1991), 2337-2342.

\bibitem{xu}
Xu Y., ``Hahn, Jacobi, and Krawtchouk polynomials of several variables,"
 Journal of Approximation Theory {\bf 195} (2015) 19-42,
 {\tt arXiv:1309.1510[math.CA]}.
 
 \bibitem{zheda}
Zhedanov A., 
``9j-symbols of the oscillator algebra and Krawtchouk polynomials in two variables,"
 J. Phys. A: Math. Gen. {\bf 30} (1997)  8337-8353.
\end{thebibliography}
\end{document}